\documentclass[12pt,reqno]{amsart}
\usepackage{cite}
\usepackage[colorlinks,citecolor=blue,urlcolor=blue,hypertexnames=false]{hyperref}
\usepackage{amssymb,amsmath,amsfonts,latexsym, ulem}
\usepackage{bm}
\usepackage[us,12hr]{datetime}
\usepackage{mathtools, mathrsfs}
\allowdisplaybreaks
\usepackage{enumerate}
\usepackage{fancyhdr}
\usepackage{appendix}

\usepackage{array,graphics,xcolor}
\numberwithin{equation}{section}
\newtheorem{dfn}{Definition}[section]
\newtheorem{thm}[dfn]{Theorem}
\newtheorem{lma}[dfn]{Lemma}
\newtheorem{hypo}{Hypothesis}
\newtheorem{conv}{Convention}
\newtheorem{ppsn}[dfn]{Proposition}
\newtheorem{cor}[dfn]{Corollary}

\newtheorem{rmrk}[dfn]{Remark}
\newcommand\at[3]{#1\,#2\Big|_{#3}}

\newcommand{\R}{\mathbb{R}}
\newcommand{\N}{\mathbb{N}}
\newcommand{\C}{\mathbb{C}}

\newcommand{\Scal}{\mathcal{S}}
\newcommand{\Hcal}{\mathcal{H}}
\newcommand{\Bcal}{\mathcal{B}}
\newcommand{\Rcal}{\mathcal{R}}
\newcommand{\Wcal}{\mathcal{W}}
\newcommand{\bh}{\mathcal{B}(\Hcal)}
\newcommand{\Tr}{\operatorname{Tr}}

\def\Im{{\mathrm{Im}\,}}

\newcommand{\Z}{\mathbb{Z}}
\renewcommand{\S}{\mathcal{S}}
\renewcommand{\H}{\mathcal{H}}

\newcommand{\mB}{\mathcal{B}}

\newcommand{\nrm}[2]{\left\|#1\right\|_{#2}}
\newcommand{\norm}[1]{\left\|#1\right\|}

\newcommand{\mW}{\mathscr{W}}

\newcommand{\bT}{\mathbf{T}}
\newcommand{\bV}{\mathbf{V}}

\newcommand{\mF}{\mathfrak{F}}

\newcommand{\abs}[1]{\left|#1\right|}

\DeclareMathOperator{\dom}{dom}
\DeclareMathOperator{\ran}{ran}

\newcommand{\nn}{\nonumber}
\begin{document}
	\title[SSF and BMV conjecture for relatively bounded perturbations]{BMV conjecture and positivity of spectral shift for relatively Schatten perturbations}
	
	\author[Pradhan]{Chandan Pradhan}
	\address{C. P., Department of Mathematics and Statistics, MSC01 1115, University of New Mexico, Albuquerque, NM 87131, USA}
	\email{chandan.pradhan2108@gmail.com, cpradhan@unm.edu}
	
	\author[Skripka]{Anna Skripka}
	\address{A. S., Department of Mathematics and Statistics, MSC01 1115, University of New Mexico, Albuquerque, NM 87131, USA}
	\email{askripka@unm.edu}
	
	
	
	\begin{abstract}
		We establish BMV-type results and prove the positivity-preserving property and uniqueness of support-normalized spectral shift functions of arbitrary order for semibounded operators with relatively Schatten–von Neumann perturbations. We also derive correlation inequalities in the relative Hilbert--Schmidt case. These results build on a new positivity-preserving property for the trace of a higher-order operator derivative along a relatively Schatten direction and are obtained using methods of multilinear operator integration and approximation.
	\end{abstract}
	
	\maketitle

	\section{Introduction}\label{sec:intro}
	
	In this paper, we establish several positivity results which, despite their apparent differences, admit a unified treatment and hold for symmetric, relatively Schatten-class perturbations of operators with essential spectrum. In particular, we prove BMV-type results and correlation inequalities, establish positivity properties of spectral shift functions, and show that the traces of certain higher-order operator derivatives are positive. Analogues of our BMV-type results and correlation inequalities for operators with discrete spectrum were obtained in \cite{DNN} (see also \cite{LNR}), while prototypes of our BMV-type results, as well as of the positivity results for spectral shift functions and traces of operator derivatives, were established in \cite{PS25} in the setting of Schatten-class perturbations of operators with essential spectrum. In the semibounded case, the weaker assumptions imposed here allow us to cover a broader class of models than those considered in \cite{PS25}, while the inclusion of operators with continuous spectrum extends the scope beyond the discrete-spectrum setting of \cite{DNN,LNR}; see Section~\ref{sec:exs} for examples.
	
	\medskip
	
	\paragraph{\bf BMV-type results}
	
	The much-studied BMV conjecture asserts that the partition function of a quantum statistical system described by a semibounded self-adjoint operator $H$ and a nonnegative perturbation $V$, regarded as a function of the coupling constant $t$, is the Laplace transform of a positive measure (see \cite[(I.3)]{BMV}). Following the resolution of the conjecture in the finite-dimensional setting \cite{Stahl}, attention shifted to the infinite-dimensional setting. It was shown in \cite[Theorem 10]{DNN} that the BMV conjecture holds for a family of self-adjoint operators $H+tV$, defined on a common dense domain in a separable complex Hilbert space $\Hcal$, under the assumption that $e^{-(H+tV)}$ is trace class, and hence the spectrum of $H+tV$ is purely discrete, for every $t\in[-1,1]$. More precisely, there exists a nonnegative measure $\mu$ on $\R$ such that
	\begin{align*}
		\Tr(e^{-H-tV})=\int_\R e^{-ts}d\mu(s),\quad t\in[-1,1].
	\end{align*}
	
	If the operators $H+tV$ have nontrivial essential spectrum, then $\Tr(e^{-H-tV})$ is not well-defined, and a correction term needs to be subtracted, which modifies the resulting integral representation. It was shown in \cite{PS25} that if $H\ge mI$ for some $m\in\R$ and $V\ge 0$ is trace class, then $t\mapsto\Tr(e^{-H-tV}-e^{-H})$ admits a L\'evy--Khintchine representation in terms of a positive measure $\mu$ and a nonpositive constant $b$, that is,
	\begin{align}
		\label{BMV1intro}
		\Tr\big(e^{-H-tV} - e^{-H}\big) = b t + \int_{(0,\infty)} \big(e^{-ts}-1\big)\, d\mu(s),\quad t\ge 0.
	\end{align}
	Analogously, it was shown in \cite{PS25} that if $r\ge 1$ and $\lambda<m$, then there exists a positive measure $\nu$ and a nonpositive constant $c$, depending on $r$ and $\lambda$, such that
	\begin{align}
		\label{BMV2intro}
		\Tr\big((H+tV-\lambda I)^{-r}-(H-\lambda I)^{-r}\big) = c t + \int_{(0,\infty)} \big(e^{-ts}-1\big)\, d\nu(s),\quad t\ge 0.
	\end{align}
	
	In this paper, we extend the representations \eqref{BMV1intro} and \eqref{BMV2intro} to perturbations that need not be trace class and may even be unbounded.
	In Theorems \ref{thm:BMV-unbdd} and \ref{thm2:BMV-unbdd}, we prove that the representation \eqref{BMV1intro} remains valid for $V$ satisfying $V(H-i)^{-p}\in\Scal^{n/p}$ for a fixed $n\in\N$ and every $p=1,\dots,n$,
	where $\Scal^\alpha$ denotes the $\alpha$th Schatten--von Neumann ideal of compact operators on $\Hcal$. (See, e.g., \cite{Simon-tr-id} for the definitions and properties of $\Scal^\alpha$.) In Theorems \ref{thm:BMV-unbdd} and \ref{thm2:BMV-unbdd}, we also prove that \eqref{BMV2intro} holds under the same relaxed assumptions on $V$, provided that $r\ge 1$ if $n=1$ and $r\ge 2n$ if $n\ge 2$.
	
	\medskip
	
	\paragraph{\bf Correlation inequalities}
	
	Correlation inequalities are among the standard tools of statistical physics; see \cite{LNR} and the references therein. The result in \cite[Theorem 7.1]{LNR} controls the second moment of an observable $V$ in the unperturbed Gibbs state $\displaystyle\Gamma=\frac{e^{-H}}{\Tr(e^{-H})}$ by its first moments in a family of perturbed Gibbs states, together with terms involving the Gibbs expectations of the double commutator and the operator norms of the double and quadruple commutators. This inequality was sharpened in \cite[Theorem 2]{DNN}, where the bound involves only a uniform first-moment estimate and the expectation of the double commutator in the unperturbed Gibbs state. To prove this refinement, the authors established the infinite-dimensional BMV-type result in \cite[Theorem 10]{DNN}. In both \cite[Theorem 7.1]{LNR} and \cite[Theorem 2]{DNN}, $H$ is assumed to satisfy $\Tr(e^{-H})<\infty$ and hence to have purely discrete spectrum; \cite[Theorem 7.1]{LNR} further assumes that $V$ is bounded and self-adjoint and that $[[H,V],V]$ is bounded, whereas \cite[Theorem 2]{DNN} allows $V$ to be symmetric and relatively $H$-bounded with relative bound strictly less than $1$ (see Definition \ref{defn:relatively bounded}) although the right-hand side of the resulting inequality may be infinite.
	
	In Theorem \ref{thm:corr-3}, we establish a second-order correlation inequality for relatively trace-class symmetric perturbations $V$ of a lower-semibounded operator $H$ with essential spectrum and bounded $[H,V]$; namely,
	\begin{align*}
		\Tr\bigl(e^{-\frac{1}{2}H} V^2  e^{-\frac{1}{2}H}\bigr) \leq   \frac12\left(\Tr\left(V e^{-H + V}\right)-\Tr\left(V e^{-H - V}\right)\right)
+ \frac{1}{4} \Tr\big(e^{-\frac{1}{2} H} [[V, H], V] e^{-\frac{1}{2} H}\big).
	\end{align*}
	In Theorem~\ref{Thm66}, we show that if $[H,V]\in\Scal^2$, then the estimate involving the double commutator reduces to one involving a single commutator; namely
	\begin{align*}
		\Tr\bigl(e^{-\frac{1}{2}H} V^2  e^{-\frac{1}{2}H}\bigr) \leq   \frac12\left(\Tr\left(V e^{-H + V}\right)-\Tr\left(V e^{-H - V}\right)\right) + \frac{1}{4} e^{-m}\|[V,H]\|_2^2.
	\end{align*}
	The latter result is the essential-spectrum analogue of \cite[Theorem 12]{DNN}.
	
	\medskip
	
	\paragraph{\bf Positivity and uniqueness of spectral shift functions}
	
	Spectral shift theory, which originated in physics in the 1940s and has since undergone substantial development, provides an alternative to the BMV approach for encoding the effect of a perturbation. In this paper, we address the sign-definiteness of spectral shift functions; for references to other results and applications see, for example, \cite{ST19}.
	
	Let $f \in \mW_n(\R)\supset C_c^{n+1}(\R)$ (see, \eqref{functiona-class-wn}). We denote by $\Rcal_n(f,H,V)$ the $n$th-order Taylor remainder defined by
	\begin{align}\label{eq:rem}
		\Rcal_n(f,H,V)
		= f(H+V) - \sum_{k=0}^{n-1} \frac{1}{k!} \at{\frac{d^k}{dt^k}}{f(H+tV)}{t=0},
	\end{align}
	where the G\^{a}teaux derivatives of operator functions are computed in the operator norm. The existence of the operator derivatives in the setting considered in this paper follows from Theorem~\ref{thm:rel_bdd_differetiation} (see also \cite[Theorem~22]{ChNuPr26}), whereas their existence in the setting of bounded perturbations is discussed, for instance, in \cite{AzCaDoSu09,Pe06,ST19}.
	
	Let $\Omega \subset \mathbb{R}$ be an interval and let $G_\Omega$ denote the closure of the convex hull of the union of the spectra of $H+tV$, that is,
	\begin{align}
		\label{def:m-M}
		G_\Omega=\overline{{\rm c.v.h.}\big(\cup_{t\in\Omega}\sigma(H + tV)\big)}.
	\end{align}

	The following benchmark result on the existence of spectral shift functions was established in \cite{Krein}, \cite{Koplienko}, and \cite{PSS} for the cases $n=1$, $n=2$, and $n \geq 3$, respectively.	
	
	\begin{thm}\label{thm:existancessf}
		Let $n \in \N$, and let $H, V$ be self-adjoint operators in $\Hcal$ such that $V \in \Scal^n$. Then, there exists a unique $\xi_{n,H,V} \in L^1(\R)$ such that
		\begin{align}\label{traceformula_inv}
			\Tr\big(\Rcal_n(f,H,V)\big)
			= \int_\R f^{(n)}(\lambda)\,\xi_{n,H,V}(\lambda)\, d\lambda
		\end{align}
		for every $f \in C_c^{n+1}(\R)$.
		The function $\xi_{n,H,V}$ is supported in the set $G_{[0,1]}$ defined in \eqref{def:m-M}.
	\end{thm}
	
	The following sign-definiteness property of the higher-order spectral shift functions whose existence is guaranteed by Theorem~\ref{thm:existancessf} was recently established in \cite[Theorem 4.1]{PS25}, thereby resolving a problem that had remained open for over a decade.
	
	\begin{thm}\label{thm:nth_order_positivity-ssf-bdd}
		Let $n \in \N$, and let $H, V$ be self-adjoint operators in $\Hcal$ such that $V \in \Scal^n$.
		Let $\xi_{n, H, V}$ be given by Theorem \ref{thm:existancessf}. Then, the following assertions hold.
		\begin{enumerate}[(i)]
			\item If $n$ is even, then $\xi_{n, H, V} \geq 0$.
			\item If $n$ is odd and $V \ge 0$ (respectively, $V\leq 0$), then $\xi_{n, H, V} \ge 0$ (respectively, $\xi_{n, H, V}\leq 0$).
		\end{enumerate}
	\end{thm}
	
	The trace formula~\eqref{traceformula_inv} was recently derived in~\cite{ChNuPr26} for symmetric, relative Schatten-class perturbations (see Theorems \ref{thm:all_order_ssf} and \ref{thm:nth_order_ssf} for details); the corresponding spectral shift functions are locally integrable and unique up to a polynomial summand.
	
	In Theorems~\ref{thm:all_oreder_positivity-ssf} and \ref{thm:nth_order_positivity-ssf}, we establish analogues of Theorem~\ref{thm:nth_order_positivity-ssf-bdd} for a symmetric perturbation $V$ of a semi-bounded self-adjoint operator $H$ such that $V$ is relatively $H$-bounded with relative bound $0$. More specifically, in Theorem~\ref{thm:all_oreder_positivity-ssf}, we prove that if, for some fixed $n \in \N$ and $m \in \R$, $H \geq mI$ (respectively, $H \leq mI$) and $V(H-i)^{-p} \in \Scal^{n/p}$ for every $p=1,\dots,n$, then there exist unique, locally integrable spectral shift functions of all orders supported on $[\Lambda_{a, b}(m), \infty)$ (respectively, $(-\infty,\widetilde\Lambda_{a, b}(m)]$) (see, \eqref{new-bd}) such that even-order spectral shift functions are nonnegative whereas the odd-order spectral shift functions have the same sign as $V$ whenever $V$ is sign-definite.
	An analogous result is established in Theorem~\ref{thm:nth_order_positivity-ssf} for the $n$th-order spectral shift function corresponding to a perturbation satisfying $V(H-i)^{-1} \in \Scal^n$. We also obtain several improvements to technical results in \cite{ChNuPr26}.
	
	\medskip
	
	\paragraph{\bf Methods}	
	
	The foundational property underlying our results described above is the positivity-preserving property of the trace of a higher-order operator derivative established in Theorems \ref{thm:Derivative-sign-bdd}, \ref{thm:Derivative-sign-rel-bdd}, and \ref{thm:Derivative-sign-rel-bdd-new}. The proofs of the BMV-type results additionally involve the positivity-preserving property of the first-order spectral shift function and the theory of Bernstein functions.
	The proof of our second-order correlation inequality for relatively trace-class perturbations additionally involves an analogue of \cite[Theorem 7.2]{LNR} obtained in Theorem~\ref{thm:corr-1}. Our approach relies, in particular, on multilinear operator integration and multistage approximation.
	
	\medskip
	
	The organization of the paper is as follows. In Section~\ref{sec:Preliminaries} we discuss preliminaries on multilinear operator integration in the setting of relatively bounded perturbations; in Section~\ref{sec:finite} we establish the positivity-preserving property of the trace of higher-order operator derivatives; in Section~\ref{sec_pos_ssf} we prove the positivity-preserving property of the spectral shift functions for relatively Schatten-class perturbations, along with their uniqueness; in Section~\ref{sec:bmv} we prove the BMV-type results; in Section~\ref{sec:cor} we obtain the essential-spectrum variants of correlation inequalities; in Section~\ref{sec:exs} we discuss examples of models satisfying the assumptions of our results.

	\section{Multilinear Operator Integration}\label{sec:Preliminaries}
	
	We recall here the fundamental properties of multilinear operator integrals from \cite{PSS, ChNuPr26} needed for the proofs of our main results.

	Let $\N=\{1,2,\dots\}$ and $\N_0=\{0\}\cup\N$. Let $\Omega \subset \mathbb{R}$ be an interval, where $\subset$ denotes inclusion that is not necessarily strict. Let $C^0(\Omega):= C(\Omega)$, the space of all continuous functions on $\Omega$. Let $n\in\N$. Let $C^n(\Omega)$ denote the space of all $n$-times continuously differentiable functions on $\Omega$. Let $C_c^n(\Omega)$ and $C_c^\infty(\Omega)$ denote the spaces of all $n$-times continuously differentiable functions and infinitely many times continuously differentiable functions on $\R$ with compact support contained in $\Omega$, respectively.
	
	The divided difference of $f\in C^n(\Omega)$ of order $n$ is
	defined recursively as follows:
	\begin{align*}
		&f^{[0]}(\lambda)=f(\lambda),\\
		&f^{[n]}(\lambda_0,\lambda_1,\ldots,\lambda_n)=\begin{cases}
			\frac{f^{[n-1]}(\lambda_0,\lambda_1,\ldots,\lambda_{n-2},\lambda_n)-f^{[n-1]}(\lambda_0,\lambda_1,\ldots,\lambda_{n-2},\lambda_{n-1})}{\lambda_n-\lambda_{n-1}} \quad &\text{if}\quad \lambda_n\neq\lambda_{n-1},\\
			\frac{\partial}{\partial \lambda}f^{[n-1]}(\lambda_0,\lambda_1,\ldots,\lambda_{n-2},\lambda)\big|_{\lambda=\lambda_{n-1}}\quad& \text{if}\quad \lambda_n=\lambda_{n-1}.
		\end{cases}
	\end{align*}

	Let $L^p(\Omega)$ denote the standard $L^p$-space on $\Omega\subset \R^d$ with respect to the Lebesgue measure and let $\ell^p(\Z^d)$ denote the standard $\ell^p$-space on $\Z^d$. We define $L^p_{\rm{loc}}(\R)$ to be the space of measurable functions $f$ on $\R$ such that $f \in L^p(K)$ for every compact set $K \subset \R$.
	
	\smallskip
	
	For $n\in\N_0$, let $\Wcal_n(\R)$ denote the Wiener class given by
	\begin{align*}
\Wcal_n(\R)&=\{f\in C^n(\R):\; f^{(n)}\in L^\infty(\R),\, \widehat{f^{(n)}}\in L^1(\R)\},
	\end{align*}
where the Fourier transform $\hat{f}$ of a nonintegrable function $f$ is understood in the sense of tempered distributions.
	Let $$u(x)=x-i,\quad x\in\R,$$ and define
	\begin{align}\label{functiona-class-wn}
		\mW_n(\R)=\{ f\in C^n(\R):\; fu^p \in \Wcal_p(\R),\, p=0,\ldots, n\}.
	\end{align}

	\noindent An elementary Fourier analytic argument (see, e.g., \cite[Lemma 4]{ChNuPr26}) yields
	\begin{align*}
		\left\{f\in C^{ n+1}(\R):\; f^{(n)}, f^{(n+1)}\in L^2(\R)\right\}\subset\Wcal_n(\R), \text{ and } C_c^{n+1}(\R)\subset\mW_n(\R) \subset \bigcap_{p=0}^{n}\, \Wcal_p(\R).
	\end{align*}

	\noindent We have the following elementary estimate.
	\begin{lma}\label{lem:simple}
		Let $a, b \in \R$ with $a < b$, let $k \in \mathbb{N}$, and let $f \in C_c^k((a,b))$. Then,
		\[
		\|f^{(j)}\|_{\infty} \le (b-a)^{k-j} \|f^{(k)}\|_{\infty},
		\quad j = 0, \dots, k.
		\]
	\end{lma}
	
	\smallskip	
	
	The following multilinear operator integral was introduced in \cite{PSS}. We consider $\Scal^\infty=\Bcal(\Hcal)$ with the operator norm.

	\begin{dfn}\label{def:moi}
		Let $n\in\N$, let $H_0,H_1,\ldots,H_n$ be self-adjoint operators in $\Hcal$, and let $f\in \Wcal_n(\R)$.
		Let $\alpha,\alpha_k\in[1,\infty]$,  $k=1,\dots,n$, satisfy $\frac{1}{\alpha}=\frac{1}{\alpha_1}+\frac{1}{\alpha_2}+\cdots+\frac{1}{\alpha_n}$ and let $E^i_{l,m}= E_{H_i}\big(\big[\frac{l}{m},\frac{l+1}{m}\big)\big)$ for $m\in \mathbb{N}$ and $l\in \mathbb{Z}$, where $E_{H_i}$ is the spectral measure of $H_i$, $i=0,\dots,n$. Define a multilinear transformation on $\Scal^{\alpha_1}\times\cdots\times\Scal^{\alpha_n}$ by
		\begin{align}\label{a1}
			&T^{H_0,\ldots,H_n}_{f^{[n]}}(V_1,V_2,\ldots,V_n)\\
			\nonumber
			&=\lim_{m\to\infty}\lim_{N\to\infty}\sum_{|l_0|,|l_1|,\ldots,|l_n|\leq N}f^{[n]}\Big(\frac{l_0}{m},\frac{l_1}{m},\ldots,\frac{l_n}{m}\Big)E_{l_0,m}^0V_1E_{l_1,m}^1
			V_2E_{l_2,m}^2\cdots V_n E_{l_n,m}^n,
		\end{align}
		where $\lim\limits_{N\to\infty}$ is evaluated in the norm $\|\cdot\|_\alpha$ when $\alpha<\infty$ and in the strong operator topology when $\alpha=\infty$; $\lim\limits_{m\to\infty}$ is evaluated in $\|\cdot\|_\alpha$.
		The transformation $T^{H_0,\ldots,H_n}_{f^{[n]}}$ is called a multilinear operator integral with symbol $f^{[n]}$.
	\end{dfn}
	The existence of the limits in \eqref{a1} for $f\in \Wcal_n(\R)$ and the corresponding norm bound
		\begin{align*}
			\|T_{f^{[n]}}^{H_0, \ldots, H_n}\| \leq \frac{1}{n!} \|\widehat{f^{(n)}}\|_1,
		\end{align*}
are justified in \cite[Lemmas 3.5, 5.1, 5.2]{PSS}. We also note that if $V_1\in\Scal^1$, then using the cyclicity of the trace, it follows from Definition \ref{def:moi} that
	\begin{align*}
		\Tr\left( T_{f^{[1]}}^{H_0, H_0}(V_1)\right) = \Tr\left(f'(H_0)V_1\right)
	\end{align*}
	for all $f\in\Wcal_1(\R)$.
	
	The following estimate is an extension of the result of \cite[Theorem 5.3 and Remark 5.4]{PSS} for $H_0=\cdots=H_n$ to the case of distinct $H_0,\dots,H_n$, whose proof is discussed in \cite[Theorem 4.3.10]{ST19}.
	We state it for a narrower set of symbols sufficient for the purpose of this paper.
	
	\begin{thm}\label{inv-thm}
		Let $n\in\mathbb{N}$ and let $\alpha,\alpha_1,\ldots,\alpha_n\in(1,\infty)$ satisfy $\tfrac{1}{\alpha_1}+\cdots+\tfrac{1}{\alpha_n}=\tfrac{1}{\alpha}$. Let $H_0,\ldots, H_{n}$ be self-adjoint operators in $\Hcal$. Assume that $V_\ell\in\mathcal{S}^{\alpha_\ell},\, 1\leq \ell\leq n$. Then, there exists a positive constant $c_{\alpha,n}$ depending only on $n,\alpha_1,\ldots,\alpha_n$ such that
		\begin{align}\label{est}
			\big\|T^{H_0,\ldots,H_{n}}_{f^{[n]}}(V_1,V_2,\ldots,V_n)\big\|_{\alpha}\leq c_{\alpha,n}\, \|f^{(n)}\|_\infty \prod\limits_{1\leq \ell \leq n}\|V_\ell\|_{\alpha_\ell}
		\end{align}
		for every $f\in C_c^{n+1}(\R)$.
	\end{thm}
The result of Theorem \ref{inv-thm} extends to multilinear operator integrals with symbol
	$$(\lambda_0,\ldots, \lambda_{n-1})\mapsto f^{[n]}(\lambda_0,\ldots, \lambda_{n-1}, \lambda_0).$$
Applying it along with the cyclicity of the trace, H\"{o}lder's inequality, and \eqref{est} yields the bound for the trace of a multilinear operator integral given below.
	
	\begin{cor}\label{thm:tr-est}
		Let $n\in\mathbb{N}$, and let $\alpha_1,\ldots,\alpha_n\in(1,\infty)$ satisfy $\tfrac{1}{\alpha_1}+\cdots+\tfrac{1}{\alpha_n}=1$. Let $H_0, H_1, \ldots, H_{n-1}$ be self-adjoint operators in $\Hcal$. Assume that $V_\ell\in\Scal^{\alpha_\ell},\, 1\leq \ell\leq n$. Then, there exists a positive constant $c_n$ depending only on $n,\alpha_1,\ldots,\alpha_n$ such that
		\begin{align*}
			\big|\Tr(T^{H_0,H_1,\ldots, H_{n-1}, H_0}_{f^{[n]}}(V_1,\ldots,V_n))\big|&\leq c_{n}\|f^{(n)}\|_\infty\prod\limits_{1\leq \ell \leq n}\|V_\ell\|_{\alpha_\ell}
		\end{align*}
		for every $f\in C_c^{n+1}(\R)$.
	\end{cor}

	The following continuity of multilinear operator integrals can be obtained by adjusting the proof of \cite[Lemma~2.5]{PS25}.
	
	\begin{lma}\label{lem:moi-cont}
		Let $n\leq p \in \N$.
		Let $H_0, \ldots, H_n, H_{0k}, \ldots, H_{nk}$ be self-adjoint operators in $\Hcal$ such that $H_{ik} \to H_i$ in the strong resolvent sense as $k\rightarrow\infty$ for $i=0,\dots,n$ and let $V_k, V \in \Scal^p$ be self-adjoint operators such that $\|V_k-V\|_p\rightarrow 0$ as $k\rightarrow\infty$. Then, for each $f \in \Wcal_n(\R)$,
		\begin{align*}
			\big\| T_{f^{[n]}}^{H_{0k},\ldots,H_{nk}}\big( (V_k)^n \big)
			- T_{f^{[n]}}^{H_0, \ldots, H_n}\big( (V)^n \big) \big\|_{p/n} \to 0
		\end{align*}
		as $k \to \infty$.
	\end{lma}

More generally, under the same assumptions on $H_0, \ldots, H_n, H_{0k}, \ldots, H_{nk}$, the conclusion holds with $(V_k)^n$ and $(V)^n$ replaced by $(V_{1k},\ldots,V_{nk})$ and $(V_1,\ldots,V_n)$, respectively, and with $\|\cdot\|_{p/n}$ replaced by $\|\cdot\|_q$, provided that $V_{jk},V_j\in\Scal^{p_j}$, $\|V_{jk}-V_j\|_{p_j}\to0$, $1\leq p_j<\infty$, $j=1,\dots,n$, and $\sum_{j=1}^n p_j^{-1}=q^{-1}\leq1$. The operators $(V_{1k},\ldots,V_{nk})$ and $(V_1,\ldots,V_n)$ need not be self-adjoint. This extension follows from the Fourier integral representation of the divided difference, H\"older's inequality, and dominated convergence.
	
	For $a, b \in \R$, set $a \wedge b = \min\{a, b\}$ and $a \vee b = \max\{a, b\}$.
	One of our tools is the following change of variables formula derived in \cite[Lemma 9]{NuSu22}.
	
	\begin{lma}\label{lem:insert-weight}
		Let $k, n\in \N\cup \{0\}$. Let $f\in C^n(\R)$ such that $\widehat{(fu^m)^{(l)}}\in L^1(\R)$ for all $m=0,\dots,k$ and $l=0,\dots,n$. Let $H_0, \ldots, H_n$ be self-adjoint operators in $\Hcal$, and let $V_1,\dots, V_n\in\bh$. Then
		\begin{align*}
			&T_{f^{[n]}}^{H_0,\dots, H_n}(V_1,\ldots, V_n)=\, \, \sum_{r=0}^{k\wedge n}(-1)^r\sum_{\substack{ j_0\geq0,\, j_1,\ldots,j_r\geq 1,\\ j_0+\cdots+j_r=k}}\\
			&\qquad \Big(T_{(fu^{k-r})^{[n-r]}}^{H_0,\ldots, H_{n-r}} (V_1, \dots, V_{n-r}) (H_{n-r}-i)^{-j_0} V_{n-r+1} (H_{n-r+1}-i)^{-j_1}\cdots V_n (H_n-i)^{-j_r}\Big).
		\end{align*}
	\end{lma}
	
	\begin{dfn}\label{defn:relatively bounded}
		Let $H$ be self-adjoint in $\H$. A linear operator $V$ with domain $\dom(V)\subseteq\H$ is called (relatively) $H$-bounded if $\dom(H)\subset \dom(V)$ and there exist $a,b\in[0,\infty)$ such that
		\begin{align}\label{Hbound}
			\norm{V\psi}\leq a\norm{H\psi}+b\norm{\psi},\quad \psi\in\dom(H).
		\end{align}
		The infimum over such numbers $a$ is called the $H$-bound of $V$.
	\end{dfn}
	
	The following properties of relatively bounded operators will be used in the sequel.
	\begin{lma}\label{lem:KR}
		Let $H$ be self-adjoint in $\H$. We use the convention $1/0=\infty$.
		\begin{enumerate}[(a)]
			\item\label{item:properties rel bdd: norm} If $V$ is relatively $H$-bounded with $H$-bound $a$, then $V(H-i)^{-1}:\H\to\H$ is bounded and there exists $b\geq0$ such that $a\leq\|V(H-i)^{-1}\|\leq a+b$.
			\item\label{item:properties rel bdd: Kato-Rellich} If $V$ is symmetric and relatively $H$-bounded with $H$-bound $a$ then
			\begin{align*}
				H+tV\quad\text{is self-adjoint on}\quad\dom(H+tV)=\dom(H)
			\end{align*}
			for all $t\in(-\tfrac1a,\tfrac1a)$.
			\item\label{item:properties rel bdd: shift} If $V$ is symmetric and relatively $H$-bounded with $H$-bound $a$, then $V$ is relatively $(H+tV)$-bounded for all $t\in \left(-\frac{1}{a}, \frac{1}{a} \right)$ with $(H+tV)$-bound $\leq \frac{a}{1-|t|a}$. Moreover, given $d>a$, there exists a constant $b\geq 0$ such that $\|V(H+tV-i)^{-1}\|\leq\frac{a+b}{1-|t|a}$ for all $t\in \left[-\frac{1}{d}, \frac{1}{d}\right]$.
		\end{enumerate}
	\end{lma}	
	\begin{proof}
		The statements \eqref{item:properties rel bdd: norm} and
		\eqref{item:properties rel bdd: Kato-Rellich} follow directly from Definition \ref{defn:relatively bounded} and the Kato--Rellich theorem \cite[Theorem V-4.3]{Kato1995},
		and are obtained in \cite[Lemma 6(a,d)]{ChNuPr26}. The first part of \eqref{item:properties rel bdd: shift} is obtained in \cite[Lemma~6(e)]{ChNuPr26}. Below we prove the existence of $b>0$ satisfying the second part of \eqref{item:properties rel bdd: shift}.
		
		Let $d>a$ and $\epsilon=(d-a)/2$. By Definition~\ref{defn:relatively bounded}, there exists $b_\epsilon\geq 0$ such that for all $t\in[-\frac1d,\frac1d]$,
		\begin{align*}
			\|V\psi\|
			\leq{}(a+\epsilon)\|H\psi\| + b_\epsilon \|\psi\|
			\leq{}(a+\epsilon)\|(H+tV)\psi\|
			+ |t|(a+\epsilon)\|V\psi\|
			+ b_\epsilon \|\psi\|
		\end{align*}
		for every $\psi\in \dom(H)$. Now let $\psi\in\Hcal$. Applying the above estimate to $(H+tV-i)^{-1}\psi$ with $t\in[-\frac1d,\frac1d]$, we obtain
		\begin{align*}
			\|V(H+tV-i)^{-1}\psi\|
			\leq\, & \frac{(a+\epsilon)}{1-|t|(a+\epsilon)}
			\|(H+tV)(H+tV-i)^{-1}\psi\| \\
			& + \frac{b_\epsilon}{1-|t|(a+\epsilon)}
			\|(H+tV-i)^{-1}\psi\|.
		\end{align*}
		By the spectral theorem, it follows that
		\begin{align*}
			\|V(H+tV-i)^{-1}\psi\|
			\leq{} & \Big(
			\frac{(a+\epsilon)}{1-|t|(a+\epsilon)}
			+ \frac{b_\epsilon}{1-|t|(a+\epsilon)}
			\Big)\|\psi\|,\quad t\in\Big[-\frac1d, \frac1d\Big].
		\end{align*}
		We further estimate
		\begin{align*}
			\|V(H+tV-i)^{-1}\psi\|\leq{} & \frac{a}{1-|t|a}\|\psi\|
			+ \Bigg(\frac{a+\epsilon}{1-|t|(a+\epsilon)} - \frac{a}{1-|t|a}
			+ \frac{b_\epsilon}{1-|t|(a+\epsilon)}\Bigg)\|\psi\|\\
			\leq{} &\frac{a}{1-|t|a}\|\psi\|
			+ \frac{\epsilon+b_\epsilon(1-|t|a)}{(1-|t|a)(1-|t|(a+\epsilon))}\|\psi\|\\
			\leq & \frac{a}{1-|t|a}\|\psi\|+ \frac{1}{1-|t|a}\,
			\frac{\epsilon + b_\epsilon}{1-|t|(a+\epsilon)}\,\|\psi\|
			\leq\frac{a+b}{1-|t|a}\|\psi\|,
		\end{align*}
		where $b= \frac{\epsilon + b_\epsilon}{1-d^{-1}(a+\epsilon)}$.
	\end{proof}
	
The following properties of relatively bounded perturbations follow from \cite[Theorem 6.4, Lemma 6.22]{Teschl2014} (see also \cite{TeschlErrata}) and will be used in the sequel.
\begin{ppsn}\label{prop-comp-lwbd}
	Let $m\in\R$. Let $H$ be self-adjoint in $\Hcal$ and let $V$ be symmetric relatively $H$-bounded operator. Then, the following assertions hold.
	\begin{enumerate}[(i)]
		\item\label{rel-comp} If $V(H-i)^{-1}$ is compact, then the $H$-bound of $V$ is $0$.
		\item\label{lower-bd} If $H\geq mI$ and the $H$-bound of $V<1$, then
		$$H+V\geq \left(m-\max\left\{a|m|+b, \frac{b}{1-a}\right\}\right)I,$$
where $a\in[0,1)$ and $b\geq0$ satisfy \eqref{Hbound}. Moreover, if $V\ge 0$, then $H+V\geq mI$.
	\end{enumerate}
\end{ppsn}

	The following second resolvent identity will be used repeatedly throughout the proofs.
	\begin{lma}{$($\!\!\cite[Lemma 7]{ChNuPr26}\,$)$}\label{lem:sri}
		Let $H$ be self-adjoint in $\Hcal$ and let $V$ be symmetric and relatively $H$-bounded with $H$-bound $<1$. For each $z\in\C$ with $\Im(z)\neq 0$, we have
		$$(H+V-z)^{-1}-(H-z)^{-1}=-(H+V-z)^{-1}V(H-z)^{-1} .$$
	\end{lma}
	
	We now summarize the construction of multiple operator integrals for relatively bounded arguments, as presented in \cite[Definition 9, Lemma 10, Theorem 15]{ChNuPr26}.
	
	\begin{dfn}\label{dfn:moi_rel_bdd}
		Let $n\in\N$, let $H_0,\ldots,H_n$ be self-adjoint in $\H$. For each $j\in\{1,\ldots,n\}$, let $V_j$ be a relatively $H_j$-bounded operator. Let $f\in \mW_n(\R)$.  Denote $$\tilde V_{j,l}:=V_{j+1}(H_{j+1}-i)^{-1}\cdots V_l(H_l-i)^{-1}\in\mB(\H), \quad \tilde V_{j,j}=I,$$ and define the multilinear operator integral
		\begin{align*}
			\bT^{H_0,\ldots,H_n}_{f^{[n]}}(V_1,\ldots,V_n):=\sum_{r=0}^n(-1)^{n-r}\sum_{0<j_1<\cdots<j_r\leq n} T^{H_0,H_{j_1},\ldots,H_{j_r}}_{(fu^r)^{[r]}}(\tilde V_{0,j_1},\ldots,\tilde V_{j_{r-1},j_r})\tilde V_{j_r,n}\,.
		\end{align*}
		The case $r = 0$ on the right-hand side of the above identity is understood as $(-1)^nf(H_0)\tilde{V}_{0,n}$. Furthermore, for all $V_1,\ldots,V_n\in\mB(\H)$ we have		$$\bT^{H_0,\ldots,H_n}_{f^{[n]}}(V_1,\ldots,V_n)=T^{H_0,\ldots,H_n}_{f^{[n]}}(V_1,\ldots,V_n).$$
For $n=0$, we use the convention $\bT_f^H()=T_f^H()=f(H)$.
	\end{dfn}	
	
We will also use the following perturbation formula for multilinear operator integrals obtained in \cite[Theorems 21 and 50]{ChNuPr26}.
	
	\smallskip
	
	\begin{lma}\label{lem:moi_relative_pert_form}
		Let $n\in\N_0$ and $j\in\{0,\ldots,n\}$. Let $H_0,\ldots,H_{n}$ be self-adjoint in $\H$, let $V_1,\ldots,V_{n}$ be operators in $\H$ such that $V_k$ is relatively $H_k$-bounded for $k=1, \ldots, n$. Let $W$ be symmetric and relatively $H_j$-bounded with $H_j$-bound $<1$. If $j > 0$, suppose that $V_j$ is relatively $H_j + W$-bounded. Then, for all $f\in\mW_{n+1}(\R)$,
\begin{align*}
&\bT^{H_0,\ldots,H_{j-1},H_j+W,H_{j+1},\ldots,H_n}_{f^{[n]}}(V_1,\ldots,V_n)
-\bT^{H_0,\ldots,H_{j-1},H_j,H_{j+1},\ldots,H_n}_{f^{[n]}}(V_1,\ldots,V_n)\\
&\qquad\qquad\qquad\qquad=\bT_{f^{[n+1]}}^{H_0,\ldots,H_{j-1},H_j+W,H_j, H_{j+1},\ldots,H_n}(V_1,\ldots,V_j,W,V_{j+1},\ldots,V_n),
\end{align*}
where, in the case $n = 0$, the left-hand side of the identity is interpreted as $f(H_0+W) - f(H_0)$.
	\end{lma}
	
	For an operator $T$ and $k\in\N_0$, we introduce the notation
	\[(T)^k :=\begin{cases}
		\underbrace{T, \dots, T}_{k \text{ times}}, & \text{if } k \neq 0,\\[2mm]
		\varnothing, & \text{if } k = 0,\end{cases}
	\]
	where $\varnothing$ denotes the empty tuple.
	
	Below we establish a higher-order differentiability of an operator function $t \mapsto f(H + tV)$ in the norm $\|\cdot\|_1$ along a relatively Schatten direction $V$ (i.e., $V(H-i)^{-1}\in\Scal^n$). This result is a counterpart of the differentiability of $t \mapsto f(H + tV)$ in the operator norm established in \cite[Theorem~22]{ChNuPr26} for a relatively bounded $V$, and will be used in the proofs of our main results.	
	
	\begin{thm}\label{thm:rel_bdd_differetiation}
		Let $n\in \N$. Let $H$ be self-adjoint in $\Hcal$, and let $V$ be symmetric and relatively $H$-bounded such that
		\begin{align*}
			V(H-i)^{-1}\in\S^{n}.
		\end{align*} Let $ f \in \mW_n(\R)$.
		Then, the map $\R \ni t \mapsto f(H + tV)$ is $n$-times differentiable in the operator norm and
		\begin{align}\label{eq:kth-dif}
			\frac{1}{k!}
			\at{\frac{d^{k}}{dt^{k}}}{f(H + tV)}{t = s}
			= \bT_{f^{[k]}}^{(H + sV)^{k+1}}\!\big((V)^k\big)
		\end{align}
		for each $k =1,\ldots, n$ and $s\in \R$.
		Moreover, the map $\R\ni s \mapsto\frac{d^{n-1}}{dt^{n-1}}f(H + tV)\big|_{t = s}$
		is differentiable in the norm $\|\cdot\|_1$ and
		\begin{align}\label{eq:nth-dif}
			\frac{1}{n!}\at{\frac{d^{n}}{dt^{n}}}{f(H + tV)}{t = s}
			=\bT_{f^{[n]}}^{(H + sV)^{n+1}}\!\big((V)^n\big), \quad  s\in \R.
		\end{align}
	\end{thm}
	
	\begin{proof}
 By Proposition \ref{prop-comp-lwbd}\eqref{rel-comp}, the $H$-bound of $V$ is $0$. The $k$-times differentiability of
\[\R \ni t \mapsto f(H + tV)\]
in the operator norm was established in \cite[Theorem~22]{ChNuPr26}. Since $V$ is relatively $H+sV$-bounded (see Lemma \ref{lem:KR}\eqref{item:properties rel bdd: shift}) and since $V(H+sV-i)^{-1}\in\S^{n}$ (by Lemma \ref{lem:sri}) for every $s\in\R$, it is enough to prove the differentiability of $t\mapsto\frac{d^{n-1}}{dt^{n-1}}f(H+sV+tV)$ in $\|\cdot\|_1$ along with the representation \eqref{eq:nth-dif} at $0$ for every $s\in\R$ and note that
\begin{align*}
\at{\frac{d^{n-1}}{dt^{n-1}}}{f(H + sV + tV)}{t=\delta} = \at{\frac{d^{n-1}}{dt^{n-1}}}{f(H + tV)}{t=s+\delta}
\end{align*}
for $s,\delta \in \R$. Thus, in the proof below we assume that $s=0$.
		
		Fix $d>0$ and $\delta\in [-d,d] \setminus\{0\}$. By \eqref{eq:kth-dif}, Lemma \ref{lem:KR}\eqref{item:properties rel bdd: shift} and Lemma \ref{lem:moi_relative_pert_form}, we have
		\begin{align}\label{eq:difference}
			\nn&\frac{1}{\delta}\Big(\at{\frac{d^{n-1}}{dt^{n-1}}}{f(H + tV)}{t = \delta} - \at{\frac{d^{n-1}}{dt^{n-1}}}{f(H + tV)}{t = 0}\Big)\\
			\nn &= \frac{(n-1)!}{\delta} \Big(\bT_{f^{[n-1]}}^{(H + \delta V)^{n}}\!\big((V)^{n-1}\big) - \bT_{f^{[n-1]}}^{(H)^{n}}\!\big((V)^{n-1}\big)\Big)\\
			&= (n-1)!\sum_{k=0}^{n-1} \bT_{f^{[n]}}^{(H)^k, H+\delta V, H, (H + \delta V)^{n-k-1}}\!\big((V)^{n}\big).
		\end{align}
		
		Next we confirm that the term with $k=0$ in \eqref{eq:difference} approaches $\bT_{f^{[n]}}^{(H)^{n+1}}\!\big((V)^{n}\big)$ in the $\|\cdot\|_1$-norm as $\delta\to 0$, while the remaining terms can be handled similarly. By Definition \ref{dfn:moi_rel_bdd},
		\begin{align}\label{eq:dif-0}
			\bT_{f^{[n]}}^{ H+\delta V, H, (H + \delta V)^{n-1}}\!\big((V)^{n}\big)
			=
			\sum_{p=0}^n(-1)^{n-p}\sum_{0<j_1<\cdots<j_p\leq n}
			T^{H_0,H_{j_1},\ldots,H_{j_p}}_{(fu^p)^{[p]}}(\tilde V_{0,j_1},\ldots,\tilde V_{j_{p-1},j_p})\tilde V_{j_p,n},
		\end{align}
		where $H_j=H+\delta V$ for $j\neq 1$ and $H_1=H$ and $\tilde V_{j,l}:=V(H_{j+1}-i)^{-1}\cdots V(H_l-i)^{-1}$. Note that $\delta V$ is relatively $H$-bounded with $H$-bound $0$. By Lemma \ref{lem:sri} and Lemma \ref{lem:KR} \eqref{item:properties rel bdd: Kato-Rellich} and \eqref{item:properties rel bdd: shift} with $a=0$ and a suitable $b\geq0$, we have
		\begin{align}
			\label{eq-est-new-1}\|(H+\delta V-i)^{-1}-(H-i)^{-1}\|
			\leq\, & |\delta|\, \|V(H-i)^{-1}\| \to 0 \text{ as } \delta\to 0,\\
			\nn\|V(H+\delta V-i)^{-1}-V(H-i)^{-1}\|_n
			\nn\leq\, & |\delta|\, \|V(H+\delta V-i)^{-1} V(H-i)^{-1}\|_n \\
			\nn\leq\, & |\delta|\, b\,\|V(H-i)^{-1}\|_n \to 0 \text{ as } \delta\to 0.
		\end{align}
		Therefore, applying Lemma \ref{lem:moi-cont} along with H\"older's inequality on the right-hand side of \eqref{eq:dif-0} yields
		\begin{align*}
			\|\cdot\|_1-\lim_{\delta\to 0}\, \bT_{f^{[n]}}^{ H+\delta V, H, (H + \delta V)^{n-1}}\!\big((V)^{n}\big)
			=\bT_{f^{[n]}}^{(H)^{n+1}}\!\big((V)^{n}\big).
		\end{align*}
		
		By a similar argument applied to every $k$ in \eqref{eq:difference}, we obtain that the limit in \eqref{eq:difference} exists in the norm $\|\cdot\|_1$ and equals
		\begin{align*}
			\frac{1}{n!}\at{\frac{d^{n}}{dt^{n}}}{f(H + tV)}{t = 0}
			= \bT_{f^{[n]}}^{(H)^{n+1}}\!\big((V)^n\big).
		\end{align*}
	\end{proof}

	\section{ Positivity preserving property of the trace of operator derivatives}\label{sec:finite}
	
	In this section, we establish the positivity-preserving property of the trace of a higher-order operator derivative with respect to a relatively Schatten direction of differentiation.
	
	Certain formulas in this paper contain the term $(H-i)^{-l}V$, which is understood in the following sense.
	\begin{conv}\label{conv1}
Let $l\in\N$ and $p\in[1,\infty)$.
Let $H$ be self-adjoint in $\Hcal$ and let $V$ be symmetric and relatively $H$-bounded. Note that $(H-i)^{-l}V\psi=(V(H+i)^{-l})^*\psi$ for every $\psi\in\dom(V)$. Since $\dom(V)$ is dense in $\H$, the operator $(H-i)^{-l}V$ has a unique bounded extension on the whole space $\H$, which is still denoted by $(H-i)^{-l}V$. If $V(H-i)^{-l}\in\S^p$, then $(H-i)^{-l}V\in\S^p$ and $\|(H-i)^{-l}V\|_p=\|V(H+i)^{-l}\|_p =\|V(H-i)^{-l}\|_p$.

	\end{conv}
	
	 We will need the following three lemmas, the first two of which can be proved similarly to \cite[Lemma~31]{ChNuPr26} and \cite[Lemma~32]{ChNuPr26}.

	\begin{lma}\label{aprlem1}
		Let $n\in \N$ and let $F$ be either $\{1,\ldots,n\}$ or $\{1\}$. Let $H$ be self-adjoint in $\Hcal$, and let $V$ be symmetric and relatively $H$-bounded such that
		\begin{align*}
			V(H-i)^{-p}\in\S^{n/p}\quad\text{for every }\, p\in F.
		\end{align*}
		Then there exists a sequence $\{V_k\}_{k\in\mathbb{N}}$ of finite-rank self-adjoint operators of the form
		$$V_k = Q_k E_H((-k,k)) V E_H((-k,k)),$$ where $Q_k$ is a finite-rank projection commuting with $E_H((-k,k)) V E_H((-k,k))$, such that the following assertions hold.
		
		\begin{enumerate}[(i)]
			\item\label{item:aprlem-0}If $F=\{1\}$, then for every $k \in\N$, $\|V_k- E_H((-k,k)) V E_H((-k,k))\|_n< \frac{1}{k^2}$.
			
			\item\label{item:aprlem-1} Given $d>0$, there exists $b\geq 0$ such that for every $t\in [-d, d]$,
			\[\|V_k(H-i)^{-1}\|,\ \|V_k(H+tV_k-i)^{-1}\|\leq b.\]
			
			\item\label{item:aprlem-2} For every $p\in F$,
			\[
			\lim_{k\to\infty} \|V_k(H-i)^{-p}-V(H-i)^{-p}\|_{n/p}=0.
			\]
			\item\label{item:aprlem-3} For every $k\in\N$ and $p\in F$,
			\[
			\|V_k(H-i)^{-p}\|_{n/p}\leq \|V(H-i)^{-p}\|_{n/p}.
			\]
		\end{enumerate}
If $V\geq0$ (respectively, $V\leq0$), the construction gives $V_k\geq0$ (respectively, $V_k\leq0$).
	\end{lma}
	
	\begin{lma}\label{aprlem2}
		Let $n, H, V$ satisfy the hypothesis of Lemma \ref{aprlem1} and let $F$ be either $\{1,\ldots,n\}$ or $\{1\}$. Let $W\in\{V, V_k, V-V_k\}$ and $R\in\{V, V_k\}$, where $\{V_k\}_{k\in\N}$ is given by Lemma \ref{aprlem1}.
		Denote
		\[H_t:=H+tR,\quad t\in \R,\]
		and
		\begin{align}\label{kk3}
			\alpha_{n,V,H}:=\max_{p\in F}\|V(H-i)^{-p}\|_{n/p}.
		\end{align}
		Then, the following assertions hold.
\begin{enumerate}[(i)]
\item\label{esta} $(H_t-i)^{-p}W$, $W(H_t-i)^{-p}\in \Scal^{n/p}$  for $p\in F$.
\item\label{estb} Given $d>0$, there exists $b\geq 0$ such that for every $t\in [-d, d]$ and $p\in F$,
\[\|(H_t-i)^{-p}W\|_{n/p}=\|W(H_t-i)^{-p}\|_{n/p}\leq 2^p\,(1\vee d)^{p}\, (1+bd)\,\alpha_{n,W,H}\,\max\limits_{\ell=0,\, p-1}\{\alpha_{n,V,H}^{\ell}\}.\]
\end{enumerate}
	\end{lma}
	
		\begin{lma}\label{lem:uniform-rel-bd}
			Let $m \in \R$. Let $H \geq mI$ (respectively, $H\leq mI$) and $V$ satisfy the hypothesis of Lemma \ref{aprlem1}.
			Suppose that $a \in [0, 1)$ and $b \in [0, \infty)$ satisfy \eqref{Hbound}, and let $W \in \{V, V_k\}$, where $\{V_k\}_{k \in \mathbb{N}}$ is given by Lemma \ref{aprlem1}.
Define
\begin{align}\label{new-bd}
\Lambda_{a, b}(m) &:= \begin{cases} m &\text{if } V\ge 0,\\
m - \max\left\{ a|m| + b, \frac{b}{1-a} \right\} &\text{otherwise},\end{cases}\\ \nonumber
\widetilde\Lambda_{a, b}(m) &:= \begin{cases} m &\text{if } V\le 0,\\
 m + \max\left\{ a|m| + b, \frac{b}{1-a} \right\} &\text{otherwise}.\end{cases}
			\end{align}
Then,
			\begin{align*}
				H + tW \geq \Lambda_{a, b}(m) I\quad (\text{respectively, } H + tW \leq \widetilde\Lambda_{a, b}(m) I) \quad \text{for all } t \in [0, 1].
			\end{align*}
		\end{lma}	
		\begin{proof}
First assume that $H\geq mI$. By Proposition \ref{prop-comp-lwbd}\eqref{rel-comp}, the $H$-bound of $tW$ is $0$ for every $t\in [0, 1]$. For every $h\in \dom(H)$, using \eqref{Hbound}, we obtain
			\begin{align*}
				\|tW h\| \leq a |t| \|H h\| + b |t| \|h\|.
			\end{align*}
			Then by Proposition \ref{prop-comp-lwbd}\eqref{lower-bd}, for each $t \in [0, 1]$, we have
			\begin{align*}
				H + tW \geq \begin{cases} mI &\text{if } W\ge 0\\ \left( m - \max\left\{ |t| (a |m| + b), \frac{|t|b}{1 - |t|a} \right\}\right) I &\text{otherwise}\end{cases} \; \geq \Lambda_{a, b}(m) I.
			\end{align*}
Now assume that $H\leq mI$. Applying the established lower bound $-H-tW\ge \Lambda_{a, b}(-m) I$ yields $H + tW \leq \widetilde\Lambda_{a, b}(m) I$ for every $t\in [0, 1]$.
		\end{proof}

	Following arguments similar to those used in the proofs of \cite[Theorems 3.3 and 3.4]{PS25}, we obtain the following positivity-preserving property, which will be used in the proof of our correlation inequality.
	
	\begin{thm}\label{thm:Derivative-sign-bdd}
		Let $n\in\N$, $H, V$ be self-adjoint operators in $\Hcal$ such that $V\in\Scal^n$. Let $\Omega\subset \R$ be an interval. Let $g\in\cap_{j=1}^n\Wcal_j(\R)$ and consider the function $\psi:\Omega\to \C$ defined by
		\begin{align*}
			\psi(s)=\Tr\!\left(\frac{d^n}{dt^n} g(H + tV)\Big|_{t=s}\right).
		\end{align*}
		If $g^{(n)}|_{G_\Omega}\ge 0$, where $G_\Omega$ is given by \eqref{def:m-M}, then the following assertions hold.
		\begin{enumerate}[(i)]
			\item If $n$ is even, then $\psi \ge 0$.
			
			\item If $n$ is odd, then $\psi \ge 0$ (respectively, $\psi \le 0$) provided $V \ge 0$ (respectively, $V \le 0$).
		\end{enumerate}
	\end{thm}
	
	The next result extends Theorem \ref{thm:Derivative-sign-bdd} to the setting of relatively Schatten-class perturbations, which will be used in the proof of the respective BMV conjecture.
	
	\begin{thm}\label{thm:Derivative-sign-rel-bdd}
		Let $m\in\R$, $n\in\N$. Let $H$ be self-adjoint in $\Hcal$ such that $H\geq mI$. Let $V$ be symmetric and relatively $H$-bounded such that $V\geq 0$ and
		\begin{align*}
			V(H-i)^{-1}\in\S^{n}.
		\end{align*}
		Let $g\in \mW_n(\R)$ and consider the function $\psi:[0,\infty)\to\C$ defined by
		\begin{align*}
			\psi(s)=\Tr\!\left(\frac{d^n}{dt^n} g(H + tV)\Big|_{t=s}\right).
		\end{align*}
		If $g^{(n)}|_{[m,\infty)}\ge 0$, then $\psi\ge 0$.
	\end{thm}
	\begin{proof}
		Let $g\in \mW_n(\R)$ be such that $g^{(n)}|_{[m,\infty)}\ge 0$. Let $\{V_k\}_{k\in\N}$ be the sequence of nonnegative finite-rank operators given by Lemma \ref{aprlem1}. Let $W\in\{V, V_k\}$. By Theorem \ref{thm:rel_bdd_differetiation}, $\frac{d^n}{dt^n} g(H + tW)\big|_{t=s}$ exists in the trace-class norm for all $s\in [0,\infty)$, and
		\begin{align}\label{eq:der-1}
			\psi(s)= n!\Tr\!\left(\bT_{g^{[n]}}^{(H+sV)^{n+1}}((V)^{n})\right)
		\end{align}
		and
		\begin{align}\label{eq:der-2}
			\psi_k(s):=\Tr\!\left(\frac{d^n}{dt^n} g(H + tV_k)\Big|_{t=s}\right)
			= n!\Tr\!\left(\bT_{g^{[n]}}^{(H+sV_k)^{n+1}}((V_k)^{n})\right).
		\end{align}
		
		Next, we show that for each fixed $s$, the sequence $\psi_k(s)$ converges to $\psi(s)$. To see this, according to Definition \ref{dfn:moi_rel_bdd}, we note that
		\begin{align}\label{eq:der-3}
			\nn
			&\Tr\!\Big(\bT_{g^{[n]}}^{(H+sV_k)^{n+1}}((V_k)^{n})\Big)\\
			&=\sum_{p=0}^n\,\sum_{\substack{j_1+\cdots+j_{p+1}= n\\j_1,\ldots, j_p{\ge 1};\, j_{p+1}\geq 0}}\, (-1)^{n-p}\, \Tr\,\Big( T^{(H+sV_k)^{p+1}}_{(gu^p)^{[p]}}(\tilde V^{j_1}_k,\ldots,\tilde V^{j_p}_k)\tilde V^{j_{p+1}}_k\Big)
		\end{align}
		where $\tilde V_k:= V_k(H+sV_k-i)^{-1}$. Let $\tilde V=V(H+sV-i)^{-1}$.
		For each $s\in[0,\infty)$, by Lemmas \ref{lem:sri}, \ref{aprlem1} and \ref{aprlem2},
		\begin{align}\label{eq-new-0}
			&\nn\|(H+sV-i)^{-1}-(H+sV_k-i)^{-1}\|_n  = |s|\, \|(H+sV_k-i)^{-1}(V-V_k)(H+sV-i)^{-1}\|_n\\
			& \leq\, |s|\,\|(V-V_k)(H+sV-i)^{-1}\|_n \longrightarrow 0 \text{ as } k\to\infty
		\end{align}
		and
		\begin{align}\label{eq-new-1}
			\nn&\|V(H+sV-i)^{-1}-V_k(H+sV_k-i)^{-1}\|_n\\
			\nn& \leq \|(V-V_k)(H+sV-i)^{-1}\|_n+  |s|\, \|V_k(H+sV_k-i)^{-1}(V-V_k)(H+sV-i)^{-1}\|_n\\
			& \leq\, \big(1+|s|\,\|V_k(H+sV_k-i)^{-1}\|\big)\,\|(V-V_k)(H+sV-i)^{-1}\|_n \longrightarrow 0 \text{ as } k\to\infty.
		\end{align}
		Therefore, by Lemma \ref{lem:moi-cont},
		\begin{align*}
			\|T^{(H+sV_k)^{ p+1}}_{(gu^p)^{[p]}}(\tilde V^{j_1}_k,\ldots,\tilde V^{j_p}_k)\tilde V^{j_{p+1}}_k-T^{(H+sV)^{ p+1}}_{(gu^p)^{[p]}}(\tilde V^{j_1},\ldots,\tilde V^{j_p})\tilde V^{j_{p+1}}\|_1 \longrightarrow 0 \text{ as } k\to \infty
		\end{align*}
		for $p=0,\ldots, n$.
		Combining the latter with \eqref{eq:der-1}, \eqref{eq:der-2}, and \eqref{eq:der-3} gives
		\begin{align}\label{eq:final_limit}
			\lim_{k\to\infty} \psi_k(s) = \psi(s).
		\end{align}
By \cite[Theorem 3.4]{PS25}, $\psi_k(s)\geq 0$ for all $k\in\N$ and $s\in [0,\infty)$. Therefore, thanks to \eqref{eq:final_limit}, we conclude that $\psi(s)\geq 0$ for all $s\in [0,\infty)$.
	\end{proof}
	
	\smallskip

	Let $n\in\N$, $m,\lambda,r\in\R$, $r\ge 2n$, $\lambda<m$. On $[m,\infty)$, consider the functions
	$$f_1(x)=-e^{-x},\quad f_2(x)=-(x-\lambda)^{-r}.$$
	Clearly, $f_i\in C^\infty([m,\infty))$ and $(f_iu^k)^{(k)}, (f_iu^k)^{(k+1)}\in L^2([m,\infty))$ for all $k\in\N_0$, $i=1,2$. We also note that $(f_iu^p)^{(k)}, (f_iu^p)^{(k+1)} \in L^2([m, \infty))$ for all $p=0,\dots,n$ and $k=0,\dots,n-1$.
	Let $g_i\in C^{\infty}(\R)$, $i=1,2$, be such that
	\begin{align}\label{two functions}
		g_i(x)=\begin{cases}
			0 & \text{ if } x\leq \frac{m+\lambda}{2}\\
			f_i(x) & \text{ if } x \geq m
		\end{cases}
	\end{align}
	for $i=1, 2$. Then, $g_i\in \cap_{j=1}^\infty \mW_j(\R)$, $i=1,2$.

	\begin{thm}\label{thm:Derivative-sign-rel-bdd-new}
		Let $n\in\N$, $m,\lambda,r\in\R$,  $r\ge 2n$, $\lambda<m$. Let $H$ be self-adjoint in $\Hcal$ such that $H\geq mI$. Let $V$ be symmetric and relatively $H$-bounded such that $V\geq 0$ and
		\begin{align*}
			V(H-i)^{-p}\in\S^{n/p} \quad\text{for }\; p=1,2,\ldots, n.
		\end{align*}
		Let $g_1, g_2$ be given by \eqref{two functions} and let $\psi_i:[0,\infty)\to\bh$ be defined by
		\begin{align*}
			\psi_i(s) = g_i(H + sV)-g_i(H), \quad i=1,2.
		\end{align*}
Then, $\psi_i(s)\in\Scal^1$ for all $s\in[0, \infty)$ and $i=1, 2$. Moreover, $\psi_i$ is infinitely many times differentiable on $[0,\infty)$ in the norm $\|\cdot\|_1$ and
		\begin{align}
			\label{eq:kth-dif-new}
			&\psi_i^{(k)}(s)= k!\, \bT_{g_i^{[k]}}^{(H + sV)^{k+1}}\!\big((V)^k\big),\\
			\label{eq-sign}
			&(-1)^{k-1}\Tr\!\big(\psi_i^{(k)}(s)\big) = (-1)^{k-1} \frac{d^k}{dt^k} \Tr\!\big(g_i(H + tV)-g_i(H)\big)\Big|_{t=s}\geq 0
		\end{align}
		for all $k\in \N$, $s\in[0,\infty)$, $i=1,2$.
	\end{thm}
	
	\begin{proof}
By Proposition \ref{prop-comp-lwbd}\eqref{rel-comp}, the $H$-bound of $V$ is $0$.  Successive applications of Lemma \ref{lem:moi_relative_pert_form}, Definition \ref{dfn:moi_rel_bdd}, Lemma \ref{aprlem2}\eqref{esta}, and the remark following Definition \ref{def:moi} yield
\begin{align*}
	\psi_i(s)& = sT_{{(g_iu)}^{[1]}}^{H+sV, H}(V(H-i)^{-1}) -  s g_i(H+sV) V(H-i)^{-1}\\
	& = s T_{{(g_iu^n)}^{[1]}}^{H+sV, H}(V(H-i)^{-n}) - s \sum_{k=0}^{n-1} (g_iu^k)(H+sV) V(H-i)^{-k-1}\\
	&= s T_{{(g_iu^n)}^{[1]}}^{H+sV, H}(V(H-i)^{-n}) - s \sum_{k=0}^{n-1} (g_iu^{k+n})(H+sV)\, ((H+sV-i)^{-n}V)\, (H-i)^{-k-1} \in \Scal^1
\end{align*}
for all $s\in[0, \infty)$ and $i=1, 2$. Now, we establish \eqref{eq:kth-dif-new} and \eqref{eq-sign} in the following two cases.
		
		\noindent{\it Case 1:} $1 \le k < n$.
		
		\noindent
		We establish \eqref{eq:kth-dif-new} by induction on $k$. The base case of the induction ($k=1$) and the induction step can be handled by the same argument, so we provide details only for the induction step. Let $l \in \N$ with $l < n - 1$, and assume that for every $k \in \{1, \ldots, l\}$, the derivative $\psi_i^{(k)}(s)$ exists in the norm $\|\cdot\|_1$ and satisfies \eqref{eq:kth-dif-new}.
		
		Let $s \in [0, \infty)$ be fixed and let $h\in\R\setminus\{0\}$ satisfy $s+h\in[0,\infty)$. By the induction hypothesis, telescoping, and Lemma~\ref{lem:moi_relative_pert_form},
		\begin{align*}
			\frac{1}{h}\Big(\psi_i^{(l)}(s+h)- \psi_i^{(l)}(s)\Big)
			\nn &= \frac{l!}{h} \Big(\bT_{g_i^{[l]}}^{(H + (s+h) V)^{l+1}}\!\big((V)^{l}\big) - \bT_{g_i^{[l]}}^{(H+sV)^{l+1}}\!\big((V)^{l}\big)\Big)\\
			&= l!\sum_{p=0}^{l} \bT_{g_i^{[l+1]}}^{(H+sV)^p,\, H+(s+h) V,\, H+sV,\, (H + (s+h) V)^{l-p}}\!\big((V)^{l+1}\big).
		\end{align*}
		To conclude the existence of the $(l+1)$st derivative of $\psi_i$ and the representation \eqref{eq:kth-dif-new}, it suffices to show that for each $p \in \{0, \ldots, l\}$, $\bT_{g_i^{[l+1]}}^{(H+sV)^p,\, H+(s+h) V,\, H+sV,\, (H + (s+h) V)^{l-p}}\!\big((V)^{l+1}\big) \in \Scal^1$
		and that it converges to $\bT_{g_i^{[l+1]}}^{(H + s V)^{l+2}}\!\big((V)^{l+1}\big)$
		in $\|\cdot\|_1$ as $h \to 0$.
		Below we include a detailed proof for the term with $p = 0$, while the remaining cases can be handled similarly.
		
		By Definition \ref{dfn:moi_rel_bdd},
		\begin{align}\label{eq-est-1}
			&\bT_{g_i^{[l+1]}}^{H+(s+h) V,\, H+sV,\, (H + (s+h) V)^{l}}\!\big((V)^{l+1}\big)\\ \nn
			=&\, (-1)^{l+1} g_i(H_0)\,\tilde V_{0,l+1} + \sum_{d=1}^{l+1}(-1)^{l+1-d}\sum_{0<l_1<\cdots<l_d\leq l+1}
			T^{H_0,H_{l_1},\ldots,H_{l_d}}_{(g_iu^d)^{[d]}}(\tilde V_{0,l_1},\ldots,\tilde V_{l_{d-1},l_d})\tilde V_{l_d,l+1}
		\end{align}
		where $H_0= H+(s+h) V$, $H_1= H+sV$, and $H_j=H+(s+h) V$ for $j>1$ and
$$\tilde V_{j,d}:=V(H_{j+1}-i)^{-1}\cdots V(H_d-i)^{-1}.$$
		Denote $l_0:=0$,
		\begin{align*}
			&\bV_{j_0,d-r}:=(H_{l_{d-r}}-i)^{-j_0}\tilde V_{l_{d-r},l_{d-r+1}}(H_{l_{d-r+1}}-i)^{-j_1}\cdots \tilde V_{l_{d-1},l_{d}}(H_{l_{d}}-i)^{-j_r}\tilde V_{l_d,l+1},\\
			&\bV_{n+j_0}:=(H_0-i)^{-n-j_0}\, V\, (H_1-i)^{-1} \tilde V_{1,l_{1}}(H_{l_{1}}-i)^{-j_1}\cdots \tilde V_{l_{d-1},l_{d}}(H_{l_{d}}-i)^{-j_d}\tilde V_{l_d,l+1},\\
			&\bV_{d-r}:=\tilde V_{l_{d-r},l_{d-r+1}}(H_{l_{d-r+1}}-i)^{-j_1}\cdots \tilde V_{l_{d-1},l_{d}}(H_{l_{d}}-i)^{-j_r}\tilde V_{l_d,l+1}.
		\end{align*}
		Applying Lemma~\ref{lem:insert-weight} to each term on the right-hand side of \eqref{eq-est-1} and extracting the term with $r=d$ from the sum $\sum_{r=0}^{(n-l-1)\wedge d}$ results in
		\begin{align}\label{eq-dif-new}
			&\bT_{g_i^{[l+1]}}^{H+(s+h) V, H+sV, (H + (s+h) V)^{l}}\!\big((V)^{l+1}\big)\\
			\nn=&\, (-1)^{l+1} g_i(H_0) \tilde V_{0,l+1} + \sum_{d=1}^{l+1} (-1)^{l+1-d}
			\sum_{0 < l_1 < \cdots < l_d \leq l+1}\sum_{r=0}^{(n-l-1)\wedge d}(-1)^r\sum_{\substack{ j_0\geq0,\, j_1,\ldots,j_r\geq 1,\\ j_0+\ldots+j_r=n-l-1}}\\
			\nn &\qquad\qquad\qquad\qquad \Big(T^{H_0,\, H_{l_1},\ldots, H_{l_{d-r}}}_{(g_iu^{n-l+d-r-1})^{[d-r]}} \bigl(\tilde V_{0,l_1}, \ldots, \tilde V_{l_{d-r-1},l_{d-r}}\bigr)\bV_{j_0,d-r}\Big)\\
			\nn=&\, (-1)^{l+1} (g_iu^{n-l})(H_0) \big((H_0-i)^{-n+l}V\big)\, (H_1-i)^{-1}\, \tilde V_{1,l+1}\\
			\nn & + \sum_{d=1}^{(n-l-1) \wedge (l+1)} (-1)^{l+1}
			\sum_{0 < l_1 < \cdots < l_d \leq l+1}\,\sum_{\substack{ j_0\geq0,\, j_1,\ldots,j_d\geq 1,\\ j_0+\ldots+j_d=n-l-1}}(g_iu^{2n-l-1})(H_0)\bV_{n+j_0}\\
			\nn& + \sum_{d=1}^{l+1} (-1)^{l+1-d}
			\sum_{0 < l_1 < \cdots < l_d \leq l+1}\,\,\sum_{ r\in\{0, \dots,(n-l-1)\wedge d\}\setminus\{d\}}(-1)^r\,\sum_{\substack{ j_0\geq0,\, j_1,\ldots,j_r\geq 1,\\ j_0+\ldots+j_r=n-l-1}}\\
			\nn &\qquad\qquad\qquad\qquad\Big( T^{H_0,\, H_{l_1},\ldots, H_{l_{d-r}}}_{(g_iu^{n-l+d-r-1})^{[d-r]}} \bigl(\tilde V_{0,l_1}, \ldots, \tilde V_{l_{d-r-1},l_{d-r}} (H_{l_{d-r}}-i)^{-j_0}\bigr)\bV_{d-r} \Big).
		\end{align}
		Observe that every term on the right-hand side of \eqref{eq-dif-new} belongs to $\Scal^1$.
		By \eqref{eq-est-new-1},
		\begin{align*}
			\|(H + (s+h)V - i)^{-1} - (H + sV - i)^{-1}\|_n \to 0\quad\text{as }h \to 0.
		\end{align*}
		For $p\in\{1,\ldots, n\}$, by using a telescoping sum, Lemma \ref{lem:sri}, and H\"older's inequality, we obtain
		\begin{align*}
			&\|V(H + (s+h)V - i)^{-p} - V(H + sV - i)^{-p}\|_{n/p}\\
			=\, & |h| \Big\|\sum_{k=0}^{p-1} V(H+(s+h)V-i)^{-k-1} \, V (H+sV-i)^{-p+k} \Big\|_{n/p}\\
			\leq\,& |h| \|V(H+(s+h)V-i)^{-1}\| \, \|V (H+sV-i)^{-p}\|_{n/p}\\
			&\qquad + |h| \sum_{k=1}^{p-1} \|V(H+(s+h)V-i)^{-k}\|_{n/k} \, \|V (H+sV-i)^{-p+k}\|_{n/(p-k)}.
		\end{align*}
		By Lemma \ref{lem:KR}\eqref{item:properties rel bdd: shift} and Lemma \ref{aprlem2}\eqref{estb}, each term in the latter sum converges to $0$ as $h\to 0$. Similarly, for every $p=1, \ldots, n$,
		\[\|(H+(s+h)V-i)^{-p}V-(H+sV-i)^{-p}V\|_{n/p} \to 0\quad\text{as }h\to 0.\]
		Applying \cite[Theorem VIII.20(b)]{ReedSimonI}, Lemma \ref{lem:moi-cont}, and H\"older's inequality on the right-hand side of \eqref{eq-dif-new} yields
		\begin{align}
			\|\cdot\|_1-\lim_{h\to 0} \bT_{g_i^{[l+1]}}^{H+(s+h) V, H+sV, (H + (s+h) V)^{l}}\!\big((V)^{l+1}\big) = \bT_{g_i^{[l+1]}}^{(H + sV)^{l+2}}\!\big((V)^{l+1}\big).
		\end{align}
		Therefore, $\psi_i^{(l+1)}(s)$ exists in $\|\cdot\|_1$ and is given by \eqref{eq:kth-dif-new} with $k=l+1$ for every $l=0,1,\dots,n-2$.
		
		By continuity of the trace, for every $k=0,\dots,n-1$, $s\in[0,\infty)$, and $i=1,2$,
		\begin{align*}
			\Tr\!\left(\psi_i^{(k)}(s)\right) = \frac{d^k}{dt^k} \Tr\!\left(g_i(H + tV)-g_i(H)\right)\Big|_{t=s}.
		\end{align*}
		
		Next we prove \eqref{eq-sign} for $1\leq k<n$. Let $\{V_l\}$ be the sequence of nonnegative finite-rank self-adjoint operators given by Lemma \ref{aprlem1}. For each $l\in\N$, define
		\begin{align*}
			\psi_{i,l}(s)= g_i(H+sV_l)-g_i(H),\quad s \in [0,\infty),\; i=1, 2.
		\end{align*}
		For $s \in [0, \infty)$ and $i \in \{1, 2\}$, by Theorem \ref{thm:rel_bdd_differetiation}
		\begin{align}\label{eq-l-k-est}
			\psi_{i,l}^{(k)}(s) = k!\, \bT_{g_i^{[k]}}^{(H + sV_l)^{k+1}}\!\big((V_l)^k\big).
		\end{align}
		Note that the transformation on the right-hand side of \eqref{eq-l-k-est}, after division by $k!$, coincides with the transformation on the left-hand side of \eqref{eq-dif-new} with $l+1=k$, $h=0$, and $V=V_l$.
		By \eqref{eq-new-0},
		\begin{align*}
			&\|(H+sV-i)^{-1}-(H+sV_l-i)^{-1}\|_n \longrightarrow 0 \quad\text{as } l\to\infty.
		\end{align*}
		For $p\in\{1,\ldots, n\}$, by telescoping, Lemma \ref{lem:sri}, and H\"older's inequality,
		\begin{align*}
			&\|V(H+sV-i)^{-p}-V_l(H+sV_l-i)^{-p}\|_{n/p} \\
			\leq\, &\|V((H+sV-i)^{-p}-(H+sV_l-i)^{-p})\|_{n/p} + \|(V-V_l)\, (H+sV_l-i)^{-p}\|_{n/p}\\
			=\, & |s|  \Big\|\sum_{k=0}^{p-1} V(H+sV-i)^{-k-1} (V-V_l)\, (H+sV_l-i)^{-p+k} \Big\|_{n/p} + \|(V-V_l)\, (H+sV_l-i)^{-p}\|_{n/p} \\
			\leq \, & |s| \, \|V(H+sV-i)^{-1}\|\,  \|(V-V_l)\,(H+sV_l-i)^{-p}\|_{n/p} + \|(V-V_l)\, (H+sV_l-i)^{-p}\|_{n/p}\\
			& + |s|\, \sum_{k=1}^{p-1} \|V(H+sV-i)^{-k}\|_{n/k}\,  \|(V-V_l)\,(H+sV_l-i)^{-p+k}\|_{n/(p-k)}.
		\end{align*}
		By Lemma \ref{lem:KR}\eqref{item:properties rel bdd: shift}, Lemma \ref{aprlem2}, and Lemma \ref{aprlem1} each term in the latter sum converges to $0$ as $l\to\infty$. Similarly, for every $p=1,\ldots, n$,
		\begin{align*}
			\|(H+sV-i)^{-p}V-(H+sV_l-i)^{-p}V_l\|_{ n/p} \longrightarrow 0\quad\text{as }
			l\to \infty.
		\end{align*}
		Applying \cite[Theorem VIII.20(b)]{ReedSimonI}, Lemma \ref{lem:moi-cont}, and H\"older's inequality on the right-hand side of \eqref{eq-l-k-est} implies that for each $s \in [0, \infty)$,
		\begin{align}\label{eq:limit-new}
			\|\cdot\|_1-\lim_{l\to \infty}\, \psi_{i,l}^{(k)}(s) = \psi_{i}^{(k)}(s), \quad i=1, 2.
		\end{align}
		By Theorem \ref{thm:Derivative-sign-bdd},
		$(-1)^{k-1}\Tr\big(\psi_{i, l}^{(k)}(s)\big)\geq 0$ for all $l\in\N$, $s\in [0, \infty)$, and $i=1, 2$. Using \eqref{eq:limit-new}, we obtain $(-1)^{k-1}\Tr\big(\psi_{i}^{(k)}(s)\big)\geq 0$ for all $s\in [0, \infty)$ and $i=1, 2$.
		
		\smallskip
		
		{\it Case 2:} $k\geq n$.
		
		\noindent
		This case is a direct consequence of Theorems \ref{thm:rel_bdd_differetiation} and \ref{thm:Derivative-sign-rel-bdd} as $V(H-i)^{-1}\in\Scal^k$ for every natural $k\ge n$.
	\end{proof}

	\section{Sign-definiteness of spectral shift functions}
	\label{sec_pos_ssf}
	In this section, we establish results on the sign-definiteness and uniqueness of spectral shift functions for semibounded $H$ and sign-definite $V$ under each of the following two hypotheses.

	\begin{hypo}\label{mainhypo-1}
		Let $n\in \N$, $n\geq 2$. Let $H$ be self-adjoint in $\Hcal$, and let $V$ be symmetric and relatively $H$-bounded such that
		\begin{align*}
			V(H-i)^{-p}\in\S^{n/p}\quad \text{for }\; p=1,\ldots,n.
		\end{align*}
	\end{hypo}

 We note that the assumption $n\geq2$ in Hypothesis \ref{mainhypo-1} is not restrictive since $V(H-i)^{-1}\in\S^1$ (Hypothesis \ref{mainhypo-1} for $n=1$) implies $V(H-i)^{-1}\in\S^2$ and $V(H-i)^{-2}\in\S^2$ (Hypothesis \ref{mainhypo-1} for $n=2$).
	\begin{hypo}\label{mainhypo-2}
		Let $n\in \N$. Let $H$ be self-adjoint in $\Hcal$, and let $V$ be symmetric and relatively $H$-bounded such that
		\begin{align*}
			V(H-i)^{-1}\in\S^{n}.
		\end{align*}
	\end{hypo}
	
	The trace formulas for relatively Schatten-class perturbations were established for the following classes of functions.
	\begin{dfn}
		Let $n\in\N$. Let $\mF_n(\R)$ denote the space of all functions $f \in C^n(\R)$ such that
		\begin{enumerate}[(i)]
			\item $\widehat{f^{(k)}u^k} \in L^1(\R)$,  $k=0,\ldots,n$,
			\item $f^{(k)}u^{k-1} \in L^1(\R)$,  $k=1,\ldots,n$.
		\end{enumerate}
	\end{dfn}	
	
	\begin{dfn}\label{def:function class Q_n}
		Let $n,k\in \mathbb{N}$. Let $\mF_n^k(\R)$ denote the space of all functions $f\in C^k(\R)$ such that
		\begin{enumerate}
			\item $fu^{2n}\in C_b(\R)$,
			\item  $f^{(l)}u^{n + l + 1}\in C_0(\R)$, $ l=1,\dots,k$,~\text{and}
			\item  $\widehat{f^{(k)}u^{k\vee n}}\in L^1(\R)$.
		\end{enumerate}
	\end{dfn}
	
	The existence of the spectral shift functions and the corresponding trace formulas stated in the next two theorems are due to~\cite{ChNuPr26}; see also Remark \ref{rem:convergence} on the proof of existence.
	
	\begin{thm}{$($\!\!\cite[Theorem 39]{ChNuPr26}\,$)$}\label{thm:all_order_ssf}
		Assume Hypothesis \ref{mainhypo-1} and the notation of \eqref{def:m-M}. Let $a\in[0,1)$ and $b\in [0,\infty)$ satisfy \eqref{Hbound}.
		Then, for each $k\in \N$, there exists a locally integrable real-valued function $\eta_k:=\eta_{k,H,V}$ such that
		\begin{align*}
			\int_\R \frac{|\eta_{k}(\lambda)|}{(1+|\lambda|)^{ n+k+\epsilon}}d\lambda\leq C_{n,k,\epsilon}\,\frac{1+b}{1-a}\,\max\limits_{\ell=1, \,k\vee n}\{\alpha_{n,V,H}^{\ell}\}
		\end{align*}	
		for all $\epsilon\in(0,1]$, where $C_{n,k,\epsilon}$ are some positive constants depending only on $n,k, \epsilon$, and $\alpha_{n,V,H}$ is given by \eqref{kk3}, and, moreover,
		\begin{align}\label{cond}
			\Tr\left(\Rcal_k(f, H, V)\right)=\int_{\R}f^{(k)}(\lambda)\eta_{k}(\lambda)d\lambda
		\end{align}
		for each $f\in\mF_n^k(\R) $. The locally integrable function $\eta_{k}$ is determined by \eqref{cond} uniquely up to a polynomial summand of degree at most $k-1$.
	\end{thm}
	
	\begin{thm}{$($\!\!\cite[Theorem 40]{ChNuPr26}\,$)$}\label{thm:nth_order_ssf}
		Assume Hypothesis \ref{mainhypo-2} and the notation of \eqref{def:m-M}. Let $a\in[0,1)$ and $b\in [0,\infty)$ satisfy \eqref{Hbound}. Then, there exists $C_n>0$ and a real-valued function $\eta_n:=\eta_{n, H, V}$
		such that
		\begin{align*}
			\int_\R |\eta_n(x)|\,\frac{dx}{(1+|x|)^{n+\epsilon}}\leq C_n\,(1+\epsilon^{-1})\, \frac{1+b}{1-a}\,\nrm{V(H-i)^{-1}}{n}^n\quad\text{for all } \epsilon >0
		\end{align*}	
		and
		\begin{align}
			\label{tff}
			\Tr\big(\Rcal_n(f, H, V)\big)=\int_\R f^{(n)}(x)\eta_n(x)\,dx\,
		\end{align}
		for every $f\in \mF_n (\R)$. The locally integrable function $\eta_n$ is determined by \eqref{tff} uniquely up to a polynomial summand of degree at most $n-1$.
	\end{thm}
	
	Let $m\in\N_0$ and $d>0$. Let $\mathcal{P}_m$ denote the space of polynomials on $\mathbb{R}$ of degree at most $m$, with the same notation used for the space of their restrictions to $(-d,d)$. For $h\in L^1((-d,d))$, let $[h]$ denote its equivalence class in the quotient space $L^1((-d,d))/\mathcal{P}_m$, given by
	$[h]:=h+\mathcal{P}_m=\{h+p:p\in\mathcal{P}_m\}$.
	We will need the following lemma.
	
	\begin{lma}\label{lem:duality}
		Let $m\in\N_0$. Let $d>0$ and $\eta\in L^1((-d,d))$. Denote $[\eta]= \eta + \mathcal{P}_m \in L^1((-d,d))/\mathcal{P}_m$. Then,
		\begin{align}\label{eq:duality-formula}
			\|[\eta]\|_{L^1((-d,d))/\mathcal P_m}
			= \sup_{\substack{ f\in C_c^\infty((-d,d))\\ \|f^{(m+1)}\|_{L^\infty((-d,d))}\leq 1}} \abs{\int_{-d}^{d}\eta(x)f^{(m+1)}(x)\,dx}.
		\end{align}
	\end{lma}
	
	\begin{proof}
		Regard $\mathcal{P}_m$ as a closed subspace of $L^1((-d,d))$.
		Its annihilator $\mathcal{P}_m^\perp$ in the space $L^\infty((-d,d))=(L^1((-d,d)))^*$ is given by
		$$\mathcal P_m^\perp=\left\{q\in L^\infty((-d,d)):\int_{-d}^d x^jq(x)\,dx=0,\ j=0,\ldots,m\right\}.$$
		By \cite[Theorem 4.9(b)]{Rudin_Func},
		\begin{align}\label{eq:moment-duality}
			&\|[\eta]\|_{L^1((-d,d))/\mathcal P_m}
			=\sup_{\substack{q\in\mathcal P_m^\perp\\\|q\|_{L^\infty((-d,d))}\leq1}}
			\abs{\int_{-d}^d \eta(x)q(x)\,dx}.
		\end{align}
		If $f\in C_c^\infty((-d,d))$, then $\int_{-d}^{d}x^jf^{(m+1)}(x)\,dx=0$ for $j=0,\ldots,m$, so
		\begin{align}
			\label{eq:first-inequality}
			\|[\eta]\|_{L^1((-d,d))/\mathcal P_m}\ge
			\sup_{\substack{f\in C_c^\infty((-d,d))\\\|f^{(m+1)}\|_\infty\leq 1}}
			\abs{\int_{-d}^d \eta(x)f^{(m+1)}(x)\,dx}.
		\end{align}
		Below we prove the reverse inequality in \eqref{eq:first-inequality}.
		
		Let $q\in\mathcal P_m^\perp$ with $\|q\|_{L^\infty((-d,d))}\leq1$ and extend $q$ and $\eta$ by $0$ on $\R\setminus(-d,d)$. Choose an even nonnegative mollifier
		$\rho\in C_c^\infty((-1,1))$ with $\int_\R\rho(x)\,dx=1$. The family $\rho_\varepsilon(x)=\varepsilon^{-1}\rho(x/\varepsilon)$, $\varepsilon>0$, is an approximate identity (see, e.g., \cite[Definition~1.2.15 and Example~1.2.17]{Grafakos2014}).
		For $N\in\N,N\geq2$, set
		\[r_N=1-1/N,\quad \varepsilon_N=d/(2N),\quad h_N=\rho_{\varepsilon_N}*\big(q(\,\cdot/r_N)\big).\]
		By the support and smoothing properties of convolution (see, e.g., \cite[Theorem~2.3.20 and Exercise~1.2.1]{Grafakos2014}),
		$h_N\in C_c^\infty((-d,d))$, since $r_Nd+\varepsilon_N<d$. By the binomial formula and the vanishing moments of $q$,
		\[
		\int_{-d}^dx^jh_N(x)\,dx = r_N\int_\R\rho_{\varepsilon_N}(z)\int_{-d}^d(r_Ny+z)^jq(y)\,dy\,dz=0,\quad j=0,\ldots,m.
		\]
		Consequently, $h_N\in\mathcal P_m^\perp$. By Young's inequality (see, e.g., \cite[Theorem~1.2.12]{Grafakos2014}), $\|h_N\|_\infty\leq1$. Note also that $h_N=f_N^{(m+1)}$, where $f_N(x)=\frac1{m!}\int_{-\infty}^x(x-y)^mh_N(y)\,dy\in C_c^\infty((-d,d))$. By the approximate identity theorem
		\cite[Theorem~1.2.19(1)]{Grafakos2014} and continuity of
		normalized dilations in $L^1(\R)$, $r_N(\rho_{\varepsilon_N}*\eta)(r_N\,\cdot)\to\eta$
		in $L^1(\R)$. Hence, by Fubini's theorem and the evenness of $\rho$,
		\[\int_{-d}^d\eta(x)h_N(x)\,dx=r_N\int_{-d}^d q(y)(\rho_{\varepsilon_N}*\eta)(r_Ny)\,dy
		\longrightarrow\int_{-d}^d\eta(y)q(y)\,dy.\]
		Consequently,
		\[\left|\int_{-d}^d\eta(x)q(x)\,dx\right|=
		\lim_{N\to\infty}\left|\int_{-d}^d\eta(x)f_N^{(m+1)}(x)\,dx\right|
		\leq\sup_{\substack{f\in C_c^\infty((-d,d))\\
				\|f^{(m+1)}\|_\infty\leq1}}\left|\int_{-d}^d\eta(x)f^{(m+1)}(x)\,dx\right|.\]
		Taking the supremum over $q\in\mathcal P_m^\perp$ with
		$\|q\|_\infty\leq1$ and applying \eqref{eq:moment-duality} gives
		\[\|[\eta]\|_{L^1((-d,d))/\mathcal P_m}=
		\sup_{\substack{q\in\mathcal P_m^\perp\\\|q\|_\infty\leq1}}
		\left|\int_{-d}^d\eta(x)q(x)\,dx\right|\leq\sup_{\substack{f\in C_c^\infty((-d,d))\\
				\|f^{(m+1)}\|_\infty\leq1}}\left|\int_{-d}^d\eta(x)f^{(m+1)}(x)\,dx\right|.\]
		Combining the latter with \eqref{eq:first-inequality} completes the proof.
	\end{proof}

	Below we establish the positivity-preserving property and uniqueness of a support-normalized representative of $\eta_{k,H,V}$ given by Theorem \ref{thm:all_order_ssf} for a semibounded operator $H$.
	
	\begin{thm}\label{thm:all_oreder_positivity-ssf}
		Let $m \in \R$. Let $H \geq mI$ (respectively, $H \leq mI$) and $V$ satisfy Hypothesis \ref{mainhypo-1}. For each $k \in \N$, let $a, b$, and $\eta_{k,H,V}$ be as in Theorem \ref{thm:all_order_ssf} and assume the notation of \eqref{kk3} and \eqref{new-bd}. Then, there exists a unique locally integrable function $\hat\eta_{k,H,V}$ supported on $ [\Lambda_{a, b}(m), \infty)$ (respectively, on $(-\infty, \widetilde\Lambda_{a, b}(m)]$) and satisfying \eqref{cond} such that $\hat\eta_{k,H,V}-\eta_{k,H,V} \in \mathcal{P}_{k-1}$. Moreover, the following assertions hold.
\begin{enumerate}[(i)]
\item If $k$ is even, then $\hat\eta_{k,H,V} \geq 0$.
\item  Let $k$ be odd. If $V \ge 0$, then $\hat\eta_{k,H,V} \ge 0$;
if $V\leq 0$, then $\hat\eta_{k,H,V}\leq 0$.
\end{enumerate}
	\end{thm}

	\begin{proof}
 We prove the result only for $H\geq mI$. The proof in the case $H \leq mI$ is completely analogous. First we prove the existence and uniqueness of the specified $\hat\eta_{k,H,V}$ for arbitrary $k$ and $V$. We then prove (i) and (ii). In the sequel, we denote positive constants depending only on $a$, $b$, $k$, $n$, $d$, $H$, $V$ by $c^{(j)}_{a,b,k,n,d,H,V}$.

 By Lemma \ref{lem:uniform-rel-bd}, $H+tV\geq \Lambda_{a, b}(m) I$ for every $t\in [0, 1]$.	It follows from \eqref{cond} that
		\begin{align*}
			\int_{\R}f^{(k)}(\lambda)\eta_{k, H, V}(\lambda)\,d\lambda = \Tr\left(\Rcal_k(f, H, V)\right) = 0
		\end{align*}
		for all $f\in C_c^\infty((-\infty, \Lambda_{a, b}(m)))$ and $k\in\N$. This implies that the $k$-th distributional derivative of $\eta_{k,H,V}$ vanishes on $(-\infty, \Lambda_{a, b}(m))$ (see, e.g., \cite[Definition~3.5, part~2, p.~32]{Grubb2009}). Hence, by \cite[Theorem~4.19]{Grubb2009}, $\eta_{k,H,V}$ coincides with a polynomial of degree at most $k-1$ almost everywhere on $(-\infty, \Lambda_{a, b}(m))$, say $P_k$. Let $\hat\eta_{k,H,V}= \eta_{k,H,V}-P_k$. Then $\hat\eta_{k,H,V}=0$ on $(-\infty, \Lambda_{a, b}(m))$ and
		\begin{align}\label{eq:new-ssf}
			\int_{\R}f^{(k)}(\lambda) \hat\eta_{k, H, V}(\lambda)\,d\lambda = \int_{\R}f^{(k)}(\lambda) \eta_{k, H, V}(\lambda)\,d\lambda
		\end{align}
		for all $f\in \mF_n^k(\R)$.
		Let $\{V_l\}_{l\in\N}$ be a sequence of finite-rank self-adjoint operators as in Lemma \ref{aprlem1}, and let $\{\eta_{k,H,V_l}:=\xi_{k,H,V_l}\}_{l\in\N}$ be the corresponding spectral shift functions given by Theorem \ref{thm:existancessf}, which satisfy $\eta_{k,H,V_l}=0$ on $(-\infty, \Lambda_{a, b}(m))$ since $H+tV_l\geq \Lambda_{a, b}(m)I$ for all $t\in [0,1]$ by Lemma \ref{lem:uniform-rel-bd}. We first show that, for all $d>0$ and $k>1$,
		\[ [\eta_{k,H,V_l}] \to [\hat\eta_{k,H,V}] \quad \text{in } L^1((-d,d))/\mathcal{P}_{k-1}.\]
		
		Fix $d>0$. Combining Lemma \ref{lem:duality} with \eqref{traceformula_inv} and \eqref{eq:new-ssf} yields
		\begin{align}\label{new-dualnorm-1}
			\nn\|[\eta_{k,H,V_l}-\hat\eta_{k,H,V}]\|_{L^1((-d,d))/\mathcal{P}_{k-1}}
			\nn&= \sup_{\substack{f\in C^{\infty}_c((-d,d))\\ \|f^{(k)}\|_\infty \le 1}}
			\left| \int f^{(k)}(\lambda)\big(\eta_{k,H,V_l} - \hat\eta_{k,H,V}\big)(\lambda)\, d\lambda \right| \\
			&= \sup_{\substack{f\in C^{\infty}_c((-d,d))\\ \|f^{(k)}\|_\infty \le 1}}
			\left| \Tr\big( \Rcal_k(f,H,V) - \Rcal_k(f,H,V_l) \big) \right|.
		\end{align}	
		 Let $f \in C^{\infty}_c((-d,d))$ and $W \in \{V, V_l\}$.  We note that $sW$ is relatively $H$-bounded with $H$-bound $0$, for every $s\in[0,1]$. Applying subsequently the integral representation for the remainder (see, e.g., the proof of \cite[Theorem 5.4.3]{ST19}) and Theorem \ref{thm:rel_bdd_differetiation} yields
		\begin{align}\label{eq:int_rep_1_rewrite}
			\Rcal_k(f,H,W)
			\nn&=\Rcal_{k-1}(f,H,W)- \frac{1}{(k-1)!}\frac{d^{k-1}}{dt^{k-1}} f(H+tW)\bigg|_{t=0}\\
			\nn&= \int_0^1 \frac{(1-s)^{k-2}}{(k-2)!}\, \left[\frac{d^{k-1}}{dt^{k-1}} f(H+tW)\bigg|_{t=s}-\frac{d^{k-1}}{dt^{k-1}} f(H+tW)\bigg|_{t=0}\right]\, ds \\
			&= (k-1)\int_0^1 (1-s)^{k-2}\, \left[\bT^{(H+sW)^{k}}_{f^{[k-1]}}((W)^{k-1})-\bT^{(H)^{k}}_{f^{[k-1]}}((W)^{k-1})\right]\, ds,
		\end{align}
		with the integral evaluated in the operator norm.
		
		\smallskip
		
		\noindent The map
		\begin{align}\label{themap}
			[0,1] \ni s \mapsto \bT^{(H+sW)^{k}}_{f^{[k-1]}}((W)^{k-1})
		\end{align}
		is strongly measurable since it arises as a derivative in the operator norm topology. Moreover, by \cite[Theorem 38]{ChNuPr26}, each value of the map \eqref{themap} lies in $\Scal^1$. Depending on whether $k \le n$ or $k > n$, by combining Lemma~\ref{aprlem1}, \cite[Lemmas~3.5, 5.1, 5.2]{PSS}, H\"older's inequality with either Definition~\ref{dfn:moi_rel_bdd} alongside \cite[Lemma~37]{ChNuPr26}, or Definition~\ref{dfn:moi_rel_bdd} alone, respectively, we conclude the uniform boundedness of the map~\eqref{themap} in $\Scal^1$. Consequently, \eqref{themap} is Bochner integrable in $\Scal^1$. Therefore by the continuity of the trace in the $\Scal^1$-norm applied in \eqref{eq:int_rep_1_rewrite}, we obtain
		\begin{align}\label{eq:int_rep_3_rewrite}
			&\frac{1}{k-1}\Big|\Tr\big(\Rcal_k(f,H,V) - \Rcal_k(f,H,V_l)\big)\Big|\\ \nn
			\leq~ &\int_0^1 \Big|\Tr\Big(\big(
			\bT^{(H+sV)^{k}}_{f^{[k-1]}}((V)^{k-1})
			- \bT^{(H+sV_l)^{k}}_{f^{[k-1]}}((V_l)^{k-1})
			\big)- \big(
			\bT^{(H)^{k}}_{f^{[k-1]}}((V)^{k-1})
			- \bT^{(H)^{k}}_{f^{[k-1]}}((V_l)^{k-1})
			\big)\Big)\Big|\, ds.
		\end{align}
		Next we analyze the differences appearing inside the trace on the right-hand side of \eqref{eq:int_rep_3_rewrite}. Using Lemma \ref{lem:moi_relative_pert_form}, we write
		\begin{align*}
			&\big(
			\bT^{(H+sV)^{k}}_{f^{[k-1]}}((V)^{k-1})
			- \bT^{(H+sV_l)^{k}}_{f^{[k-1]}}((V_l)^{k-1})
			\big)- \big(
			\bT^{(H)^{k}}_{f^{[k-1]}}((V)^{k-1})
			- \bT^{(H)^{k}}_{f^{[k-1]}}((V_l)^{k-1})
			\big)\\
			& = \underbrace{\big(
				\bT^{(H+sV)^{k}}_{f^{[k-1]}}((V)^{k-1})
				- \bT^{(H+sV_l)^{k}}_{f^{[k-1]}}((V)^{k-1})
				\big)}_{:=S_1} + \underbrace{\big(
				\bT^{(H+sV_l)^{k}}_{f^{[k-1]}}((V)^{k-1})
				- \bT^{(H+sV_l)^{k}}_{f^{[k-1]}}((V_l)^{k-1})
				\big)}_{:=S_2}\\
			&- \underbrace{\big(
				\bT^{(H)^{k}}_{f^{[k-1]}}((V)^{k-1})
				- \bT^{(H)^{k}}_{f^{[k-1]}}((V_l)^{k-1})
				\big)}_{:=S_3}\\
			=& S_1  + S_2 - S_3.
		\end{align*}
		Observe that
		\begin{align*}
			S_2 - S_3= \sum_{r=0}^{k-2}\Big(
			\bT^{(H+sV_l)^{k}}_{f^{[k-1]}}((V_l)^r, V-V_l, (V)^{k-r-2})- \bT^{(H)^{k}}_{f^{[k-1]}}((V_l)^r, V-V_l, (V)^{k-r-2})\Big).
		\end{align*}
		Since the summands in $S_2 - S_3$ are of the same type as $S_1$, we provide a detailed  estimation of only $\Tr(S_1)$. We decompose $S_1=s(T_1+T_2+T_3)$, where
		\begin{align*}
			T_1 &= \bT^{H+sV,\,H+sV_l,\,(H+sV)^{k-1}}_{f^{[k]}}(V-V_l, (V)^{k-1}), \\
			T_2 &= \bT^{H+sV_l,\,(H+sV)^{k-1},\,H+sV_l}_{f^{[k]}}((V)^{k-1}, V-V_l), \\
			T_3 &= \sum_{r=0}^{k-3}
			\bT^{(H+sV_l)^{r+1},\,H+sV,\,H+sV_l,\,(H+sV)^{k-r-3},\,H+sV_l}_{f^{[k]}}
			((V)^{r+1}, V-V_l, (V)^{k-r-2}).
		\end{align*}
		Note that the terms $T_2$ and $T_3$ are of the same nature as $T_1$, in the sense that the same self-adjoint operator appears in the first and last positions of the superscript of a multilinear operator integral defining each of them. Thus, we provide a detailed estimation of only $\Tr(T_1)$.
		
		By Definition \ref{dfn:moi_rel_bdd}, we have
		\begin{align}\label{eq:T_1-0}
			T_1\nn &= \bT^{H+sV,\, H+sV_l,\,(H+sV)^{k-1}}_{f^{[k]}}(V - V_l,\,( V)^{k-1}) \\
			&= (-1)^kf(H_0) \tilde V_{0,k}+ \sum_{p=1}^{k} (-1)^{k-p}
			\sum_{0 < l_1 < \cdots < l_p \leq k}
			T^{H_{0},\, H_{l_1},\cdots, H_{l_{p}}}_{(f u^p)^{[p]}}
			\bigl(\tilde V_{0,l_1}, \ldots, \tilde V_{l_{p-1},l_p}\bigr)\tilde V_{l_p,k},
		\end{align}
		where  $H_0=H+sV$, $H_1=H+sV_l,\; H_j=H+sV$ for $j\neq 1$, $V_1=V-V_l,\; V_j=V$ for $j= 2,\dots,k$, and $\tilde V_{j,l}:=V_{j+1}(H_{j+1}-i)^{-1}\cdots V_l(H_l-i)^{-1}.$
		\medskip
		
		{\it Case 1: $1<k<n$.}
		
		\noindent
		Applying Lemma~\ref{lem:insert-weight}  to each term in \eqref{eq:T_1-0} yields
		\begin{align}\label{eq:T_1-1}
			T_1
			\nn &= (-1)^k f(H_0) \tilde V_{0,k} + \sum_{p=1}^{k} (-1)^{k-p}
			\sum_{0 < l_1 < \cdots < l_p \leq k}\,\,\sum_{r=0}^{(n-k)\wedge p}(-1)^r\,\sum_{\substack{ j_0\geq0,\, j_1,\ldots,j_r\geq 1,\\ j_0+\ldots+j_r=n-k}}\\
			\nn &\times \Big( T^{H_0,\, H_{l_1},\ldots, H_{l_{p-r}}}_{(fu^{n-k+p-r})^{[p-r]}} \bigl(\tilde V_{0,l_1}, \ldots, \tilde V_{l_{p-r-1},l_{p-r}}\bigr) (H_{l_{p-r}}-i)^{-j_0}\\
			&\qquad\qquad\qquad \tilde V_{l_{p-r},l_{p-r+1}}(H_{l_{p-r+1}}-i)^{-j_1}\cdots \tilde V_{l_{p-1},l_{p}}(H_{l_{p}}-i)^{-j_r}\tilde V_{l_p,k}\Big).
		\end{align}
		Below we estimate four cases in \eqref{eq:T_1-1}.
		Let $ \alpha_{n,V-V_l,H}$ be given by \eqref{kk3}.
		
		\smallskip
		
		(a) By the functional calculus and definition of $\tilde V_{0,k}$,
		\begin{align*}
			&f(H_0) \tilde V_{0,k}\\
			=& f(H+sV) (V-V_l)(H+sV_l-i)^{-1} V(H+sV-i)^{-1}\cdots V(H+sV-i)^{-1}\\
			=&\big((fu^{n-k+1})(H+sV)\big)(H+sV-i)^{-n+k-1} (V-V_l)(H+sV_l-i)^{-1} \Big(V(H+sV-i)^{-1}\Big)^{k-1}.
		\end{align*}
		Therefore, it follows from Hypothesis \ref{mainhypo-1}, Lemma \ref{aprlem2} together with H\"older's inequality and Lemma \ref{lem:simple} that
		$f(H_0) \tilde V_{0,k}\in \Scal^1$ and
		\begin{align}
			\nn &\|f(H_0) \tilde V_{0,k}\|_1\\
			\nn\leq &\|fu^{n-k+1}\|_\infty\, \|(H+sV-i)^{-n+k-1} (V-V_l)\|_{\frac{n}{n-k+1}}\|V(H+sV-i)^{-1}\|_n^{k-1}\\
			\nn\leq&\|f^{(k)}\|_\infty\,(2d)^k\, (1+d^2)^{(n-k+1)/2}\|(H+sV-i)^{-n+k-1} (V-V_l)\|_{\frac{n}{n-k+1}}\, \|V(H+sV-i)^{-1}\|_n^{k-1}\\
			\nn\leq&\|f^{(k)}\|_\infty\, c^{(0)}_{a,b, k, n, d, H, V}\, \alpha_{n,V-V_l,H}.
		\end{align}

		(b) Let $(p,r)$ be such that $p\neq r$ and $(r,l_p)\neq (0,k)$. Denote
		$$\bV= \tilde V_{l_{p-r},l_{p-r+1}}(H_{l_{p-r+1}}-i)^{-j_1}\cdots \tilde V_{l_{p-1},l_{p}}(H_{l_{p}}-i)^{-j_r}\tilde V_{l_p,k}.$$
		By Theorem \ref{inv-thm} and H\"older's inequality,
		\begin{align*}
			&\Big\|T^{H_0,\, H_{l_1},\ldots, H_{l_{p-r}}}_{(fu^{n-k+p-r})^{[p-r]}} \bigl(\tilde V_{0,l_1}, \ldots, \tilde V_{l_{p-r-1},l_{p-r}}\bigr) (H_{l_{p-r}}-i)^{-j_0} \bV\Big\|_1\\
			&=\Big\|T^{H_0,\, H_{l_1},\ldots, H_{l_{p-r}}}_{(fu^{n-k+p-r})^{[p-r]}} \bigl(\tilde V_{0,l_1}, \ldots, \tilde V_{l_{p-r-1},l_{p-r}}(H_{l_{p-r}}-i)^{-j_0}\bigr)  \bV\Big\|_1\\
			&\leq \|(fu^{n-k+p-r})^{(p-r)}\|_\infty\,  \tilde c^{(1)}_{a, b, k, n, d, H, V}\, \alpha_{n,V-V_l,H},
		\end{align*}
		which, by the Leibniz differentiation rule combined with Lemma~\ref{lem:simple}, implies that
		\begin{align*}
			\Big\|T^{H_0,\, H_{l_1},\ldots, H_{l_{p-r}}}_{(fu^{n-k+p-r})^{[p-r]}} \bigl(\tilde V_{0,l_1}, \ldots, \tilde V_{l_{p-r-1},l_{p-r}}\bigr) (H_{l_{p-r}}-i)^{-j_0} \bV\Big\|_1\leq\|f^{(k)}\|_\infty\,  c^{(1)}_{a, b, k, n, d, H, V}\, \alpha_{n,V-V_l,H}.
		\end{align*}
		
		(c) If $r\neq p$ and $(r, l_p)=(0, k)$, then we have two subcases. If $p\geq 2$, then applying Corollary \ref{thm:tr-est}, the Leibniz differentiation rule, and Lemma~\ref{lem:simple} yields
		\begin{align*}
			&\Big|\Tr\Big(T^{H_0,\, H_{l_1},\ldots, H_{l_{ p}}}_{(fu^{n-k+p})^{[p]}} \bigl(\tilde V_{0,l_1}, \ldots, \tilde V_{l_{ p-1},k}\bigr) (H_{l_{p}}-i)^{-(n-k)}\Big)\Big|\\
			=&\Big|\Tr\Big(T^{H_0,\, H_{l_1},\ldots, H_{l_{p}}}_{(fu^{n-k+p})^{[p]}} \bigl(\tilde V_{0,l_1}, \ldots, \tilde V_{l_{ p-1},k}(H_{l_{p}}-i)^{-(n-k)}\bigr) \Big)\Big|\\
			\leq& \|(fu^{n-k+p})^{(p)}\|_\infty\,  \tilde c^{(2)}_{a, b, k, n, H, V}\, \alpha_{n,V-V_l,H}\\
			\leq& \|f^{(k)}\|_\infty\, c^{(2)}_{a, b, k, n, d, H, V}\, \alpha_{n,V-V_l,H}.
		\end{align*}
If $p=1$, then using  $\Tr\big(T_{g^{[1]}}^{H_0,H_0}(B)\big)=\Tr(g'(H_0)B)$ for $g=fu^{n-k+1}$ and $B\in\Scal^1$ (see the remark following Definition \ref{def:moi}) and applying H\"older's inequality followed by Lemma~\ref{lem:simple} yields the bound as above.
	
	(d) If $r=p$, the argument is analogous to the one in (a).
		
		\smallskip
		
		\noindent Thus, the above estimates imply that
		\begin{align*}
			|\Tr(T_1)|\leq \|f^{(k)}\|_\infty \, c^{(3)}_{a, b, k, n, d, H, V}\, \alpha_{n,V-V_l,H}.
		\end{align*}

		By similar arguments, we also obtain that
		\begin{align*}
			|\Tr(T_j)|\leq \|f^{(k)}\|_\infty \, c^{(4)}_{a, b, k, n, d, H, V}\, \alpha_{n,V-V_l,H}
		\end{align*}
		for $j=2,3$. Therefore,
		\begin{align*}
			\big|\Tr\big(S_1
			\big)\big|\leq \|f^{(k)}\|_\infty \, c^{(5)}_{a, b, k, n, d, H, V}\, \alpha_{n,V-V_l,H}.
		\end{align*}
		Similarly,
		\begin{align*}
			\big|\Tr\big(S_2-S_3\big)\big|\leq \|f^{(k)}\|_\infty \, c^{(6)}_{a, b, k, n, d, H, V}\, \alpha_{n,V-V_l,H}.
		\end{align*}
		
		{\it Case 2: $k\geq n$.}
		
		\noindent
		By Definition \ref{dfn:moi_rel_bdd}, Hypothesis \ref{mainhypo-1} and \cite[Lemmas~3.5, 5.1, 5.2]{PSS}, $S_1, S_2, S_3\in\Scal^1$.
		By arguments analogous to those in the proof of {\it Case 1},
		\begin{align*}
			\Big|\Tr\big(S_1+S_2-S_3\big)\Big| \leq \|f^{(k)}\|_\infty \,  c^{(7)}_{a, b, k, n, d, H, V}\, \alpha_{n,V-V_l,H}.
		\end{align*}
		
		Therefore, combining {\it Case 1} and {\it Case 2}, together with \eqref{new-dualnorm-1}, \eqref{eq:int_rep_3_rewrite}, and the property of $\{V_l\}_{l\in\N}$, yields
		\begin{align}\label{eq:n-1}
			\|[\eta_{k,H,V_l}-\hat\eta_{k,H,V}]\|_{L^1((-d,d))/\mathcal{P}_{k-1}} \leq c_{a, b, k, n, d, H, V}\, \alpha_{n,V-V_l,H} \longrightarrow 0
		\end{align}
		as $l\to\infty$ for all $k>1$.
\medskip
		
Let now $k=1$. Applying (56)--(59) in the proof of \cite[Theorem 36]{ChNuPr26} to the function $\eta_{1,H,V_l}-\hat\eta_{1,H,V}$ yields
\begin{align*}
		\sup_{\substack{ f\in C_c^\infty((-d,d))\\ \|f'\|_{L^\infty((-d,d))}\leq 1}} \abs{\int_{-d}^{d}(\eta_{1,H,V_l}-\hat\eta_{1,H,V})f'(x)\,dx}\rightarrow 0
\text{ as } l\rightarrow\infty.
		\end{align*}
Applying Lemma \ref{lem:duality} yields
		\begin{align}\label{eq:n-2}
			\|[\eta_{1,H,V_l}-\hat\eta_{1,H,V}]\|_{L^1((-d,d))/\mathcal{P}_0}\to 0 \text{ as } l\to \infty
		\end{align}
		for all $d>0$.
		
		Fix $d>|\Lambda_{a, b}(m)|$. Since $\|\cdot\|_{L^1((-d,d))}$ and $\|\cdot\|_{L^1((-d,\Lambda_{a, b}(m)))}$
		are norms on the finite-dimensional space $\mathcal P_{k-1}$,
		there exists $C_{k,m,d}>0$ such that
		$\|P\|_{L^1((-d,d))}\leq C_{k,m,d}\|P\|_{L^1((-d,\Lambda_{a, b}(m)))}$
		for every $P\in\mathcal P_{k-1}$.
		Since $\eta_{k,H,V_l}-\hat\eta_{k,H,V}=0$ almost everywhere on $(-d,\Lambda_{a, b}(m))$, for every
		$P\in\mathcal P_{k-1}$,
		\begin{align*}
			&\|\eta_{k,H,V_l}-\hat\eta_{k,H,V}\|_{L^1((-d,d))}\\
			&\leq \|\eta_{k,H,V_l}-\hat\eta_{k,H,V}-P\|_{L^1((-d,d))}+\|P\|_{L^1((-d,d))}\\
			&\leq \|\eta_{k,H,V_l}-\hat\eta_{k,H,V}-P\|_{L^1((-d,d))}
			+C_{k,m,d}\|\eta_{k,H,V_l}-\hat\eta_{k,H,V}-P\|_{L^1((-d,\Lambda_{a, b}(m)))}\\
			&\leq(1+C_{k,m,d})\|\eta_{k,H,V_l}-\hat\eta_{k,H,V}-P\|_{L^1((-d,d))}.
		\end{align*}
		Taking the infimum over $P\in\mathcal P_{k-1}$ and combining it with \eqref{eq:n-1} and \eqref{eq:n-2} gives
		\[\|\eta_{k,H,V_l}-\hat\eta_{k,H,V}\|_{L^1((-d,d))}
		\leq(1+C_{k,m,d})\|[\eta_{k,H,V_l}-\hat\eta_{k,H,V}]\|_{L^1((-d,d))/\mathcal P_{k-1}}
		\longrightarrow 0.\]
		Since $d>|\Lambda_{a, b}(m)|$ in the above proof was arbitrary, restriction to smaller intervals yields
		\begin{align}\label{eq:L1-conv}
			\|\eta_{k,H,V_l} - \hat\eta_{k,H,V}\|_{L^1((-d,d))} \to 0 \text{ as } l\to\infty
		\end{align}
		for all $d>0$ and $k\in\N$.

If another $\xi_{k, H, V} \in L^1_{loc}(\R)$ satisfies \eqref{cond} and is supported on $[\Lambda_{a, b}(m), \infty)$, Theorem \ref{thm:all_order_ssf} and \cite[Theorem~4.19]{Grubb2009} imply that $\xi_{k, H, V} - \hat\eta_{k, H, V} = Q$ for some $Q \in \mathcal{P}_{k-1}$. Since $\xi_{k, H, V} - \hat\eta_{k, H, V} = 0$ on $(-\infty, \Lambda_{a, b}(m))$, we obtain $Q = 0$, proving the uniqueness of $\hat\eta_{k, H, V}$.

\medskip	

 (i) For even $k$, since $\eta_{k,H,V_l}\geq 0$ by Theorem \ref{thm:nth_order_positivity-ssf-bdd}, we obtain $\hat\eta_{k,H,V}\geq 0$.
\medskip

 (ii) If $k$ is odd and $V\geq 0$ (respectively, $V\leq 0$), then by Lemma \ref{aprlem1}, $V_l\geq 0$ (respectively, $V_l\leq 0$) for all $l\in\N$. Therefore, by Theorem \ref{thm:nth_order_positivity-ssf-bdd}, $\eta_{k,H,V_l}\geq 0$ (respectively, $\eta_{k,H,V_l}\leq 0$) for all $l\in\N$; consequently, from \eqref{eq:L1-conv} we obtain $\hat\eta_{k,H,V}\geq 0$ (respectively, $\hat\eta_{k,H,V}\leq 0$).
		
\end{proof}

	By suitably adjusting the proof of Theorem \ref{thm:all_oreder_positivity-ssf},
	we obtain the sign-definiteness and uniqueness of a support-normalized representative of $\eta_{n,H,V}$ given by Theorem~\ref{thm:nth_order_ssf} for a semibounded operator $H$.
	
	\begin{thm}\label{thm:nth_order_positivity-ssf}
	Let $m \in \R$. Let $n\in\N$. Let $H \geq mI$ (respectively, $H \leq mI$) and $V$ satisfy Hypothesis \ref{mainhypo-2}. Let $a, b$, and $\eta_{n,H,V}$ be as in Theorem \ref{thm:nth_order_ssf} and assume the notation of \eqref{new-bd}. Then, there exists a unique locally integrable function $\hat\eta_{n,H,V}$ supported on
$[\Lambda_{a, b}(m), \infty)$ (respectively, on $(-\infty, \widetilde\Lambda_{a, b}(m)]$) and satisfying \eqref{tff} such that $\hat\eta_{n,H,V}-\eta_{n,H,V}\in\mathcal{P}_{n-1}$. Moreover, the following assertions hold.
\begin{enumerate}[(i)]
\item If $n$ is even, then $\hat\eta_{n,H,V} \geq 0$.
\item Let $n$ be odd. If $V \ge 0$, then $\hat\eta_{n,H,V} \ge 0$;
if $V\leq 0$, then $\hat\eta_{n,H,V}\leq 0$.
\end{enumerate}
	\end{thm}

	\begin{rmrk}\label{rem:convergence}

 While the results of \cite[Theorem~36]{ChNuPr26} and \cite[Theorem~40]{ChNuPr26} (as well as its prototype \cite[Theorem~4.1]{NuSk22}) hold, their proofs omit a technical detail, which we supply below, retaining the original notation.

Let $n$ be fixed and let $u(\lambda)=\lambda-i$. For every $p\in\{0,\dots,n-1\}$ and $k\in\N$, let $\breve\eta_{p,k}\in L^1_{loc}(\R)$ be a function satisfying
$$\Tr(\widetilde R^p_{n,H,f}(V_k))=\int(fu^p)^{(p+1)}(x)\breve\eta_{p,k}(x)\,dx$$
(see \cite[(55), (56) and Proposition 4.4]{NuSk22}). One of the intermediate steps in the proofs of \cite[Theorem~4.1]{NuSk22} and \cite[Theorems~36 and~40]{ChNuPr26} is to show that there exists $\breve{\eta}_{p}\in L^1_{loc}(\R)$ such that
\begin{align}\label{eq:fix-1}
&\lim_{k\to\infty}\int_\R (fu^p)^{(p+1)}(x)\breve{\eta}_{p,k}(x)\,dx
=\int_\R (fu^p)^{(p+1)}(x)\breve{\eta}_{p}(x)\,dx
\end{align}
for all $f\in C_c^{n+1}(\R)$. In the notation of \cite[Eq.~(60)]{ChNuPr26}, $\breve{\eta}_{0,k}=\eta_{k,1}$ and $\breve{\eta}_0=\eta_{1}$.
		
Let $a>0$. By Lemma~\ref{lem:duality} and arguments presented in the proof of \cite[Proposition~4.4]{NuSk22},
		\begin{align*}
			\|[\breve{\eta}_{p, k}] - [\breve{\eta}_{p, m}]\|_{L^1((-a,a))/\mathcal{P}_{p}}&= \sup_{\substack{ f\in C_c^\infty((-a,a))\\ \|f^{(p+1)}\|_{L^\infty((-a,a))}\leq 1}} \abs{\int_{-a}^{a}\left(\breve{\eta}_{p, k}-\breve{\eta}_{p, m}\right)f^{(p+1)}(x)\,dx}\\
			 &= \sup_{\substack{ f\in C_c^\infty((-a,a))\\ \|(fu^p)^{(p+1)}\|_{L^\infty((-a,a))}\leq 1}} \abs{\int_{-a}^{a}\left(\breve{\eta}_{p, k}-\breve{\eta}_{p, m}\right)(fu^p)^{(p+1)}(x)\,dx}\\
			&= \sup_{\substack{ f\in C_c^\infty((-a,a))\\ \|(fu^p)^{(p+1)}\|_{L^\infty((-a,a))}\leq 1}} \big| \Tr(\widetilde R_{n, H, f}^p(V_k) - \widetilde R_{n, H, f}^p(V_m))\big|
		\end{align*}
		which converges to $0$ as $k, m \to \infty$. Consequently, $\{[\breve{\eta}_{p, k}]\}_{k\in\N}$ is a Cauchy sequence in\\ $L^1((-a,a))/\mathcal{P}_{p}$ for all $a>0$. It was erroneously claimed in \cite{ChNuPr26} and \cite{NuSk22} that $\{\breve\eta_{p, k}\}_{k\in\N}$ is Cauchy in $L^1_{loc}(\R)$. Below we confirm that there exists a sequence $\{\hat\eta_{p,k}\}_{k\in\N}\subset L^1_{loc}(\R)$ such that $\hat\eta_{p, k} \in [\breve{\eta}_{p, k}]$, and $\{\hat\eta_{p, k}\}_{k\in\N}$ is Cauchy in $L^1((-a,a))$ for all $a>0$.
		
To prove the latter fact for each fixed $p$, we first justify that for each $k$ there exists a unique polynomial $Q_{p, k} \in \mathcal{P}_{p}$ satisfying
		\begin{align}
			\label{eq:moment-1}
			\int_{-1}^{1} x^l (\breve{\eta}_{p, k}(x) - Q_{p, k}(x)) \, dx = 0,
			\quad l = 0, \ldots, p.
		\end{align}
		We look for a solution to \eqref{eq:moment-1} in the form
		$Q_{p, k}(x) = \sum_{j=0}^{p} a_{jk} x^j$, which converts \eqref{eq:moment-1} to
		\begin{align}
			\label{eq:moment-2}
			\sum_{j=0}^{p} \left( \int_{-1}^{1} x^{l+j} \, dx \right) a_{jk} = \int_{-1}^{1} x^l \breve{\eta}_{p, k}(x) \, dx,
			\quad l = 0, \ldots, p.
		\end{align}
		The system \eqref{eq:moment-2} can be written in the matrix form $M a_k = b_k$, where $M=\big[\,\int_{-1}^1 x^{l+j}\,dx\,\big]_{l,j=0}^{p}$, $a_k=[a_{jk}]_{j=0}^{p}$, and $b_k=\left[\int_{-1}^{1}x^l\, \breve{\eta}_{p, k}(x)\,dx\right]_{l=0}^{p}$.
		If $c=[c_j]_{j=0}^{p}\in\C^{p+1}$, then $$\left<Mc,c\right>=\int_{-1}^1\Big|\sum_{j=0}^{p}c_jx^j\Big|^2\,dx
		=\Big\|\sum_{j=0}^{p}c_jx^j\Big\|^2.$$ By the linear independence of the monomials $1,x,\dots,x^p$ in $L^2((-1,1))$, $M$ is strictly positive-definite, so $M$ is invertible. Thus, there exists a unique solution $a_k$ to \eqref{eq:moment-2}.
		
Since $\|P\|_s:= \sum\limits_{j=0}^{p}\big|\int_{-1}^{1}x^jP(x)\,dx\big|$ defines a norm on the finite-dimensional space $\mathcal{P}_p$, there exists a constant $C_{p,a}>0$ such that
\begin{align}\label{eq:norm-eqv}
\|P\|_{L^1((-a,a))}\leq C_{p,a} \|P\|_s \, \quad \text{ for all } P\in \mathcal{P}_p.
\end{align}
		Let $\hat\eta_{p, k} :=  \breve{\eta}_{p, k} - Q_{p, k}\in[\breve{\eta}_{p, k}] \in L^1_{loc}(\R)/\mathcal{P}_p$ and $a\geq 1$. Then, for every $P\in \mathcal{P}_p$, we obtain
		\begin{align*}
			&\|\hat\eta_{p, k}- \hat\eta_{p, l}\|_{L^1((-a,a))}\\
		   =& \|(\breve{\eta}_{p, k}-\breve{\eta}_{p, l}-P) -(Q_{p, k}-Q_{p, l}-P)\|_{L^1((-a,a))}\\
		\leq& \|\breve{\eta}_{p, k}-\breve{\eta}_{p, l}-P\|_{L^1((-a,a))} + \|Q_{p, k}-Q_{p, l}-P\|_{L^1((-a,a))}\\
		\leq& \|\breve{\eta}_{p, k}-\breve{\eta}_{p, l}-P\|_{L^1((-a,a))} + C_{p, a} \, \|Q_{p, k}-Q_{p, l}-P\|_s \quad  (\text{by } \eqref{eq:norm-eqv})\\
		  = &\|\breve{\eta}_{p, k}-\breve{\eta}_{p, l}-P\|_{L^1((-a,a))} + C_{p, a}\, \sum_{j=0}^{p}
			\abs{\int_{-1}^{1}x^j (Q_{p, k}-Q_{p, l}-P)(x)\,dx}\\
		   =& \|\breve{\eta}_{p, k}-\breve{\eta}_{p, l}-P\|_{L^1((-a,a))} + C_{p, a}\, \sum_{j=0}^{p}
			\abs{\int_{-1}^{1}x^j (\breve{\eta}_{p, k}-\breve{\eta}_{p, l}-P)(x)\,dx}\quad  (\text{by } \eqref{eq:moment-1})\\
		\leq& \|\breve{\eta}_{p, k}-\breve{\eta}_{p, l}-P\|_{L^1((-a,a))} + C_{p, a}\,(1+p) \|\breve{\eta}_{p, k}-\breve{\eta}_{p, l}-P\|_{L^1((-1,1))}\\
		\leq& (1+ (1+p)C_{p, a})\, \|\breve{\eta}_{p, k}-\breve{\eta}_{p, l}-P\|_{L^1((-a,a))}.
		\end{align*}
		Taking the infimum over all $P\in\mathcal P_p$ and using the definition of the quotient norm, we obtain
		\begin{align*}
			\|\hat\eta_{p, k}- \hat\eta_{p, l}\|_{L^1((-a,a))} \leq (1+ (1+p)C_{p, a})\, \|[\breve{\eta}_{p, k}] - [\breve{\eta}_{p, l}]\|_{L^1((-a,a))/\mathcal{P}_{p}} \longrightarrow 0 \quad \text{as } k, l\to\infty.
		\end{align*}
		Since $a\geq 1$ in the above proof was arbitrary, restriction to smaller intervals yields
		\begin{align*}
			\|\hat\eta_{p, k}- \hat\eta_{p, l}\|_{L^1((-a,a))}  \longrightarrow 0 \quad \text{as } k, l\to\infty
		\end{align*}
		for all $a>0$.
		
		Thus, $\{\hat\eta_{p, k}\}_{k\in\N}$ is a Cauchy sequence in $L^1_{loc}(\R)$. Consequently, there exists a $\breve{\eta}_p\in L^1_{loc}(\R)$ such that $\hat\eta_{p, k}\to \breve\eta_p$ in $L^1_{loc}(\R)$, which further implies that
		\begin{align*}
			\lim_{k\to\infty} \int_\R (fu^p)^{(p+1)}(x) \breve\eta_{p, k}(x) \, dx = & \lim_{k\to\infty} \left(\int_\R (fu^p)^{(p+1)}(x) \hat\eta_{p, k}(x)\, dx + \int_\R (fu^p)^{(p+1)}(x) Q_{p, k}(x)\, dx\right)\\
			=&\lim_{k\to\infty} \int_\R (fu^p)^{(p+1)}(x) \hat\eta_{p, k}(x)\, dx
=\int_\R (fu^p)^{(p+1)}(x) \breve\eta_{p}(x)\, dx
		\end{align*}
		for every $f\in C_c^{n+1}(\R)$, proving \eqref{eq:fix-1}. 		
	\end{rmrk}

	\section{BMV conjecture for a relatively Schatten perturbation}\label{sec:bmv}

	\begin{thm}\label{thm:BMV-unbdd}
		Let $n\in\N$, $n\geq 2$, $m,\lambda,r\in\R$,  $r\ge 2n$, $\lambda<m$. Let $H$ be self-adjoint in $\Hcal$ such that $H\geq mI$. Let $V$ be symmetric and relatively $H$-bounded such that $V\geq 0$ and
		\begin{align*}
			V(H-i)^{-p}\in\Scal^{n/p}\quad \text{ for }\; p=1,\ldots, n.
		\end{align*}
		On $[m,\infty)$, let $f_1(x)=-e^{-x}$, $f_2(x)=-(x-\lambda)^{-r}$ and let $g_1, g_2$ be given by \eqref{two functions}.  Then, $g_i(H+tV) - g_i(H)\in\Scal^1$ and there exist $b_i \geq 0$ and positive measures $\mu_i$ on $(0,\infty)$ satisfying
		\begin{align}\label{summablitymu}
			\int_{(0,\infty)} (1 \wedge s) \, d\mu_i(s) < \infty
		\end{align}
		such that
		\begin{align}\label{fBt}
			\Tr\big(g_i(H+tV) - g_i(H)\big) = b_i\,t + \int_{(0,\infty)} \big(1-e^{-ts}\big)\, d\mu_i(s)
		\end{align}
		for every $t \in [0,\infty)$ and $i=1,2$. The pair $(b_i,\mu_i)$ is uniquely determined by the function on the left-hand side of \eqref{fBt} for $i=1,2$.
	\end{thm}
	
	\begin{proof}
By Theorem \ref{thm:Derivative-sign-rel-bdd-new}, $g_i(H+tV) - g_i(H)\in\Scal^1$ for every $t\in[0, \infty)$ and $i=1, 2$. For $i=1, 2$, consider the following functions
		\begin{align*}
			\phi_i(t)= \Tr\big(g_i(H+tV) - g_i(H)\big), \quad t\in [0,\infty).
		\end{align*}
		By Theorem \ref{thm:Derivative-sign-rel-bdd-new}, for every $k\in\N$ and $i=1,2$,
		\begin{align*}
			(-1)^{k-1}\, \phi_i^{(k)}(t) \geq 0, \quad t\in [0, \infty).
		\end{align*}
		We note that $g_1, g_2 \in \mF_n^1(\R)$. By Theorem \ref{thm:all_order_ssf} and Theorem \ref{thm:all_oreder_positivity-ssf}, there exists a nonnegative locally integrable function $\eta_{1, H,tV}$, supported on $[m, \infty)$, satisfying
		\begin{align*}
			\phi_i(t)=\int_{[m,\infty)} g_i'(x)\,\eta_{1,H,tV}(x)\,dx \geq 0,\quad i=1,2,\; t\in[0,\infty).
		\end{align*}
		Thus, $\phi_i$, $i=1,2$, are Bernstein functions on
		$[0, \infty)$. Applying \cite[Theorem 3.2]{berstein_book} to each $\phi_i$ completes the proof.
	\end{proof}
	
	\begin{thm}\label{thm2:BMV-unbdd}
		Let $m,\lambda,r\in\R$,  $r\ge 1$, $\lambda<m$. Let $H$ be self-adjoint in $\Hcal$ such that $H\geq mI$. Let $V$ be symmetric and relatively $H$-bounded such that $V\geq 0$ and $V(H-i)^{-1}\in\Scal^{1}$. On $[m,\infty)$, let $f_1(x)=-e^{-x}$, $f_2(x)=-(x-\lambda)^{-r}$ and let $g_1, g_2$ be given by \eqref{two functions}. Then $g_i(H+tV) - g_i(H)\in\Scal^1$ for all $t\in [0, \infty)$, and there exist $b_i \geq 0$ and positive measures $\mu_i$ on $(0,\infty)$ satisfying \eqref{summablitymu} such that the trace formula \eqref{fBt} holds. The pair $(b_i,\mu_i)$ is uniquely determined by the function on the left-hand side of \eqref{fBt} for $i=1,2$.
	\end{thm}
	
	\begin{proof}
		The proof is completely analogous to that of Theorem \ref{thm:BMV-unbdd}, with application of Theorems \ref{thm:Derivative-sign-rel-bdd}, ~\ref{thm:nth_order_ssf}, ~\ref{thm:nth_order_positivity-ssf} in place of Theorems \ref{thm:Derivative-sign-rel-bdd-new}, ~\ref{thm:all_order_ssf}, ~\ref{thm:all_oreder_positivity-ssf}, respectively.
	\end{proof}
	
	The following result is obtained by a suitable adjustment of the proof of \cite[Theorem~5.4]{PS25}.
	
	\begin{thm}
		\label{bmvsn}
		Let $n\in\N$, $m,\lambda,r\in\R$,  $r\ge 2n$, $\lambda<m$. Let $F$ be either $\{1,\ldots,n\}$ or $\{1\}$. Let $H$ be self-adjoint in $\Hcal$ such that $H\geq mI$. Let $V$ be symmetric and relatively $H$-bounded such that $V\geq 0$ and
		\begin{align*}
			V(H-i)^{-p}\in\S^{n/p} \text{ for } p\in F.
		\end{align*}
		Let $g_1, g_2$ be given by \eqref{two functions}. Then, for every integer $k\geq n$, and also for $1\leq k<n$ when $F=\{1,\ldots,n\}$,
		there exist finite measures $\mu_{i,k}$ on $[0,\infty)$, $i=1,2$ satisfying $(-1)^{k-1}\mu_{i, k}\ge 0$ such that
		\begin{align}\label{fBtn}
			\Tr\!\left(	\frac{d^{k}}{d\tau^{k}} g_i(H + \tau V)\Big|_{\tau=t}\right) = \int_{[0,\infty)}\, e^{-ts}\,d\mu_{i,k}(s)
		\end{align}
		for every $t \in [0,\infty)$.
	\end{thm}

	\section{Correlation inequality}\label{sec:cor}
	In this section, we establish counterparts of the discrete-spectrum correlation inequalities in \cite[Theorems~2 and~12]{DNN} for operators with essential spectrum.
	
	One of our tools is the following counterpart of \cite[Theorem~7.2(i),(ii)]{LNR} for operators with essential spectrum. In its formulation, we assume Convention~\ref{conv1} and additionally adopt the convention below for commutators. We will write $[H,V]\in\Scal^{2}$ if $H$ is a self-adjoint operator in $\H$, $V$ is an $H$-bounded symmetric operator satisfying $V(E_H([-n, n])\Hcal)\subset \dom(H)$ for all $n\in\N$, and the commutator $[H,V]$ initially defined on $\{E_H([-n,n])h: n\in\N, h\in\Hcal\}$
(or on $\dom(H)$ when $V$ is bounded) extends uniquely to an operator in $\Scal^{2}$. Let $\bh$ denote the algebra of bounded linear operators on $\Hcal$. The same domain convention applies to $[H,V]\in\bh$. For $t>0$ and $H$ bounded below, the operator $e^{-tH}V$, initially defined on $\dom(V)$, has the unique bounded extension $(Ve^{-tH})^*$, which we continue to denote by $e^{-tH}V$. Furthermore, $e^{-H/2}V^2e^{-H/2}$ denotes $(Ve^{-H/2})^*(Ve^{-H/2})$, and weighted double commutators are interpreted by expanding them into these bounded factors.	
	
	\begin{thm}\label{thm:corr-1} Let $H$ be self-adjoint in $\Hcal$ such that $H\geq mI$ for some $m\in\R$. Let $V$ be symmetric and relatively $H$-bounded such that $V(H-i)^{-1}\in\Scal^2$.
		\begin{enumerate}[(i)]
			\item The function $s \mapsto \Tr\left(e^{-\frac{1-s}{2} H}V e^{-sH} V e^{-\frac{1-s}{2} H}\right)$ is convex on $[0, 1]$ and attains its minimum at $s = 1/2$. In particular, for every $s\in[0,1]$,
			\begin{align}\label{eq:conv-ineq-1}
				0 \le \Tr\left(e^{-\frac{1}{4}H} V e^{-\frac{1}{2} H} V e^{-\frac{1}{4} H}\right) \le \Tr\left(e^{-\frac{1-s}{2} H}V e^{-sH} V e^{-\frac{1-s}{2} H}\right) \le \Tr\left(e^{-\frac{1}{2} H} V^2 e^{-\frac{1}{2} H}\right).
			\end{align}
			
			\item Assume, in addition, that either of the following two properties is satisfied:
			\begin{itemize}
				\item $[V,H]\in \Scal^2$,
				\item $[V,H]\in \bh$ and $V(H-i)^{-1}\in\Scal^1$.
			\end{itemize}
			Then, for every $s\in[0,1]$,
			\begin{align} \label{eq:conv-ineq-2}
				0 \leq\Tr\left(e^{-\frac{1}{2} H} V^2 e^{-\frac{1}{2} H}\right) - \Tr\left(e^{-\frac{1-s}{2} H}V e^{-sH} V e^{-\frac{1-s}{2} H}\right)
				\leq \frac{1}{4} \Tr\Big( e^{-\frac{1}{2} H} [[V, H], V] e^{-\frac{1}{2} H}\Big).
			\end{align}
		\end{enumerate}
	\end{thm}
	
	\begin{proof}
		Define \[\phi(s):=e^{-\frac{1-s}{2}H}Ve^{-sH}Ve^{-\frac{1-s}{2}H},\quad s\in [0,1].\] If $s > 0$, then
		\begin{align*}
			\phi(s)= e^{-\frac{1-s}{2} H} V e^{-sH} V e^{-\frac{1-s}{2} H}
			= e^{-\frac{1-s}{2} H} V(H-i)^{-1} (H-i) e^{-sH} (H-i) (H-i)^{-1} V e^{-\frac{1-s}{2} H}.
		\end{align*}
		By the spectral theorem, $(H-i)e^{-sH}(H-i) = (H-i)^2 e^{-sH} \in \bh$ for every $s \in (0,1]$. Since $(H-i)^{-1}V$, $V(H-i)^{-1} \in \Scal^2$ by the assumption, \[\phi(s)\in\Scal^1,\quad s \in (0,1].\]
		If $s=0$, then
		\begin{align*}
			\phi(0)= e^{-\frac{1}{2} H} V V e^{-\frac{1}{2} H}
			= \left(e^{-\frac{1}{2} H}(H-i)\right) \left((H-i)^{-1} V\right) \left(V(H-i)^{-1}\right) \left((H-i)e^{-\frac{1}{2} H}\right) \in \Scal^1.
		\end{align*}
		Define
		\[f(s) := \Tr\left(\phi(s)\right),\quad s \in [0, 1].\]
		The proof proceeds in four steps, gradually moving from the most restrictive assumptions to the least restrictive ones.
		\smallskip

		{\it Case 1:} $H=H^*\in\bh$ and $V=V^* \in \Scal^1$.
		
		\noindent
		Let $\{P_k\}$ be the sequence of finite-rank projections provided by \cite[Lemma 3.2]{PS25}, and define $H_k= P_k H P_k$ and $V_k= P_k V P_k$. By \cite[Lemma 3.2]{PS25}, $P_k \uparrow I$, $H_k \to H$ in the strong resolvent sense, and $V_k \to V$ in the norm $\|\cdot\|_1$.
		
		For each $s \in [0,1]$, define
		\begin{align}\label{phi-k}
			&\phi_k(s) := e^{-\frac{1-s}{2} H_k} V_k e^{-s H_k} V_k e^{-\frac{1-s}{2} H_k},\\
			\label{f-k}
			&f_k(s) := \Tr(\phi_k(s)) = \Tr\left(V_k e^{-s H_k} V_k e^{-(1-s) H_k}\right).
		\end{align}
		The second equality in \eqref{f-k} follows by the cyclicity of the trace.
		Since $H_k \to H$ in the strong resolvent sense, by \cite[Theorem VIII.20(b)]{ReedSimonI}, $e^{-s H_k} \to e^{-s H}$ strongly for each $s \in [0,1]$. Consequently, $\phi_k(s) \to \phi(s)$ in $\|\cdot\|_1$ pointwise on $[0,1]$, implying that $f_k(s) \to f(s)$ pointwise on $[0,1]$.
		
		By a straightforward calculation, for each $s \in [0,1]$,
		\begin{align*}
			&f_k'(s)=-\Tr\left([V_k,H_k]e^{-sH_k}V_ke^{-(1-s)H_k}\right)
			=\Tr\left(V_ke^{-sH_k}[V_k,H_k]e^{-(1-s)H_k}\right),\\
			&f_k''(s)= \Tr\left( e^{-s H_k} [V_k, H_k] e^{-(1-s) H_k} [H_k, V_k] \right).
		\end{align*}
		As noted in \cite[Theorem 7.2(i)]{LNR},
		\begin{align*}
			f_k''(s) = \big\|e^{-\frac{1-s}{2} H_k} [H_k, V_k] e^{-\frac{s}{2} H_k}  \big\|_2^2\geq 0
		\end{align*}
		for $s\in[0,1]$, so $f_k$ is convex on $[0,1]$ for each $k\in\N$. As noted in \cite[Theorem 7.2(ii)]{LNR},
		\begin{align} \label{eq:conv-ineq-3}
			0 \le \Tr\left(e^{-\frac{1}{4}H_k} V_k e^{-\frac{1}{2} H_k} V_k e^{-\frac{1}{4} H_k}\right)
			\le f_k(s)\le \Tr\left(e^{-\frac{1}{2} H_k} V_k^2 e^{-\frac{1}{2} H_k}\right),\quad s\in[0,1],
		\end{align}
		due to the convexity and symmetry $f_k(s)=f_k(1-s)$, as well as
		\begin{align} \label{eq:conv-ineq-4}
			0 \le \Tr\left(e^{-\frac{1}{2} H_k} V_k^2 e^{-\frac{1}{2} H_k}\right) - f_k(s)
			\le \frac{1}{4} \Tr\left(e^{-\frac{1}{2} H_k} [[V_k, H_k], V_k] e^{-\frac{1}{2} H_k}\right),\quad s\in[0,1],
		\end{align}
		due to $f_k(0)-f_k(s)\le f_k(0)-f_k(1/2)=-\int_0^{1/2} f_k'(t)dt\le -f'_k(0)/2$.
		
		Since the convexity is preserved under pointwise limits, $f$ is convex on $[0,1]$. By the symmetry $f(s)=f(1-s)$ and convexity on $[0,1]$, $f$ attains its minimum at $s = 1/2$. Passing to the limit $k \to \infty$ in \eqref{eq:conv-ineq-3} yields \eqref{eq:conv-ineq-1}. Applying H\"older's inequality and noticing that $[V_k, H_k] \to [V, H]$ in $\|\cdot\|_1$ yields the convergence of $[[V_k, H_k], V_k]$  to $[[V, H], V]$ in $\|\cdot\|_1$. Passing to the limit $k \to \infty$ in \eqref{eq:conv-ineq-4} yields \eqref{eq:conv-ineq-2}.
		
		\medskip

		{\it Case 2:} $H \geq mI$, $V=V^* \in \Scal^1$, and $[V, H] \in \bh$.
		
		\noindent
		Set $P_k = E_H((-k,k))$, $H_k = P_k H$, and $V_k = P_k V P_k$. Let $\phi_k$ and $f_k$ be defined as in \eqref{phi-k} and \eqref{f-k}, respectively. By {\it Case 1},
		$f_k$ is convex on $[0,1]$, and the pair $(H_k, V_k)$ satisfies the inequalities \eqref{eq:conv-ineq-3} and \eqref{eq:conv-ineq-4} for every $s \in [0,1]$.
		
		Since $P_k \uparrow I$, it follows that $H_k \to H$ in the strong resolvent sense and $V_k \to V$ in $\|\cdot\|_1$. Since $[V, H] \in \bh$, we obtain
		\begin{align*}
			[V_k, H_k] = P_k [V, H] P_k \to [V, H]
		\end{align*}
		strongly as $k \to \infty$. Consequently, $[[V_k, H_k], V_k]\rightarrow[[V, H], V]$ in $\|\cdot\|_1$.  An argument similar to that in {\it Case 1}  completes the proof under the assumptions of {\it Case 2}.
		
		\medskip

		{\it Case 3:} $H\geq mI$, $V$ is symmetric and $H$-bounded, and $V(H-i)^{-1}\in\Scal^2$.
		
		\noindent
		Set $W=V(H-i)^{-1}$ and $R=(H-i)^{-1}V$ and note that
		\begin{align*}
			f(s) &= \Tr\left(e^{-\frac{1-s}{2} H}V(H-i)^{-1} (H-i)e^{-sH}(H-i) (H-i)^{-1} V e^{-\frac{1-s}{2} H}\right)\\
			&= \Tr\left(e^{-\frac{1-s}{2} H}W (H-i)^2e^{-sH} R e^{-\frac{1-s}{2} H}\right),\quad s\in (0,1].
		\end{align*}

		Set $P_k=E_H((-k,k))$ and $\hat{V}_k=P_k V P_k$. By a computation similar to \cite[Eq.~(34)]{ChNuPr26}, $\hat{V}_k \in \Scal^2$ for each $k \in \N$.
		Let $\{V_k\}$ be a sequence of finite-rank self-adjoint operators provided by Lemma \ref{aprlem1} and defined by $V_k=Q_k\hat{V}_k$, where $\{Q_k\}$ is a sequence of finite-rank projections such that $Q_k\hat{V}_k=\hat{V}_kQ_k$ and $\|V_k-\hat{V}_k\|_2<\frac{1}{k^2}$. Set $W_k=V_k(H-i)^{-1}$ and $R_k=(H-i)^{-1}V_k$. By Lemma \ref{aprlem1}, the sequences $\{W_k\}$ and $\{R_k\}$ are bounded in $\bh$ and converge to $W$ and $R$, respectively, in the norm $\|\cdot\|_2$.
		
		For each $s\in[0,1]$, define
		\begin{align}\label{gkCase3}
			g_k(s):= \Tr\left( e^{-\frac{1-s}{2} H} W_k (H-i)^2e^{-sH} R_k e^{-\frac{1-s}{2} H}\right)=\Tr\left( e^{-\frac{1-s}{2} H} V_k e^{-sH} V_k e^{-\frac{1-s}{2} H}\right).
		\end{align}
		By \eqref{eq:conv-ineq-1} in {\it Case 2} applied to the pair $(H,V_k)$, each function $g_k$ is convex on $[0,1]$ and satisfies
		\begin{align*}
			0 \le \Tr\left(e^{-\frac{1}{4}H} V_k e^{-\frac{1}{2} H} V_k e^{-\frac{1}{4} H}\right)
			\le g_k(s)
			\le \Tr\left(e^{-\frac{1}{2} H} V_k^2 e^{-\frac{1}{2} H}\right),\quad s\in[0,1],
		\end{align*}
		Combining the latter with \eqref{gkCase3} and
		\begin{align*}
			g_k(0)&=\Tr\left(e^{-\frac{1}{2}H}W_k(H-i)^2 R_k e^{-\frac{1}{2}H}\right)=\Tr\left( e^{-\frac{1}{2} H} V_k V_k e^{-\frac{1}{2} H}\right)\\
			&= \Tr\left( e^{-\frac{1}{2} H}(H-i) R_k W_k (H-i)e^{-\frac{1}{2} H}\right)
		\end{align*}
		yields	
		\begin{align} \label{eq:conv-ineq-5}
			0 \le \Tr\left(e^{-\frac{1}{4}H} W_k (H-i)^2 e^{-\frac{1}{2} H} R_k e^{-\frac{1}{4} H}\right)
			\le g_k(s)
			\le \Tr\left(e^{-\frac{1}{2} H}(H-i) R_k W_k (H-i)e^{-\frac{1}{2} H}\right).
		\end{align}
		Since $e^{-\frac{1-s}{2} H}, (H-i)^2e^{-sH}\in \bh$ for each $s\in(0,1]$ and $W_k \to W$ and $R_k \to R$ in $\|\cdot\|_2$, $g_k(s)$ converges to $f(s)$ pointwise on $(0,1]$. A similar argument shows that $g_k(0)\to f(0)$. Consequently, $f$ is convex on $[0,1]$.  Passing to the limit $k\to\infty$ in \eqref{eq:conv-ineq-5} yields
		\begin{align*}
			0 \le \Tr\left(e^{-\frac{1}{4}H} W (H-i)^2 e^{-\frac{1}{2} H} R e^{-\frac{1}{4} H}\right)
			\le f(s)
			\le \Tr\left(e^{-\frac{1}{2} H}(H-i) R W(H-i) e^{-\frac{1}{2} H}\right),
		\end{align*}
		which establishes \eqref{eq:conv-ineq-1}.
		
		Next we prove \eqref{eq:conv-ineq-2} under the additional assumption $[V, H]\in \Scal^2$.
		Since $Q_k\hat{V}_k=\hat{V}_kQ_k$ and $HP_k=P_k H$,
		\begin{align}\label{eq-comm-bdd}
			[V_k, H] = [Q_k \hat{V}_k, H] = Q_k \hat{V}_k H - H \hat{V}_k Q_k = Q_k \hat{V}_k H P_k -  H P_k \hat{V}_k Q_k  \in \bh.
		\end{align}
		Therefore, by \eqref{eq:conv-ineq-2} in {\it Case 2} applied to the pair $(H,V_k)$,
		\begin{align} \label{eq:conv-ineq-6}
			0 \le \Tr\left(e^{-\frac{1}{2} H} V_k^2 e^{-\frac{1}{2} H}\right) - g_k(s)
			\le \frac{1}{4} \Tr\left(e^{-\frac{1}{2} H} [[V_k, H], V_k] e^{-\frac{1}{2} H}\right).
		\end{align}
		Next we justify that
		\begin{align}\label{eq:conv-ineq-7}
			\lim_{k\to\infty} \Tr\left(e^{-\frac{1}{2} H} [[V_k, H], V_k] e^{-\frac{1}{2} H}\right) = \Tr\left(e^{-\frac{1}{2} H} [[V, H], V] e^{-\frac{1}{2} H}\right).
		\end{align}
		The operator on the left-hand side of \eqref{eq:conv-ineq-7} can be rewritten as
		\begin{align}\label{eq:conv-ineq-8}
			\nn e^{-\frac{1}{2} H} [[V_k, H], V_k] e^{-\frac{1}{2} H}
			=&\, e^{-\frac{1}{2} H} [V_k, H] V_k e^{-\frac{1}{2} H} - e^{-\frac{1}{2} H} V_k [V_k, H] e^{-\frac{1}{2} H}\\
			\nn= &\, e^{-\frac{1}{2} H}\big( [V_k, H]- [\hat{V}_k, H]\big) V_k e^{-\frac{1}{2} H} - e^{-\frac{1}{2} H} V_k \big([V_k, H]-[\hat{V}_k, H]\big) e^{-\frac{1}{2} H}\\
			&\qquad\qquad + e^{-\frac{1}{2} H} [\hat{V}_k, H] V_k e^{-\frac{1}{2} H} - e^{-\frac{1}{2} H} V_k [\hat{V}_k, H] e^{-\frac{1}{2} H}.
		\end{align}
		Note that
		\begin{align}\label{eq-new-2}
			[\hat{V}_k, H] = P_k[V, H]P_k \xrightarrow{\|\cdot\|_2} [V, H]\quad \text{and}\quad [\hat{V}_k, H]^* = [H, \hat{V}_k] \xrightarrow{\|\cdot\|_2} [H, V].
		\end{align} Moreover, $V_k e^{-\frac{1}{2} H} = W_k (H-i)e^{-\frac{1}{2} H} \to W (H-i)e^{-\frac{1}{2} H} = V e^{-\frac{1}{2} H}$ in $\|\cdot\|_2$ and likewise $e^{-\frac{1}{2} H}V_k \to e^{-\frac{1}{2} H} V$ in $\|\cdot\|_2$. Combining these limits gives
		\begin{align*}
			&\|\cdot\|_1-\lim_{k\to\infty}\Big( e^{-\frac{1}{2} H} [\hat{V}_k, H] V_k e^{-\frac{1}{2} H} - e^{-\frac{1}{2} H} V_k [\hat{V}_k, H] e^{-\frac{1}{2} H}\Big)\\
			&=\, e^{-\frac{1}{2} H} [V, H] V e^{-\frac{1}{2} H} - e^{-\frac{1}{2} H} V [V, H] e^{-\frac{1}{2} H}.
		\end{align*}
		Proceeding similarly to \eqref{eq-comm-bdd}, we obtain
		\begin{align}\label{eq-new-3}
			\nn \|[V_k, H]- [\hat{V}_k, H]\|_2 &= \big\|\big( Q_k \hat{V}_k-\hat{V}_k \big) P_k H P_k - P_k H P_k \big( \hat{V}_k Q_k-\hat{V}_k \big)\big\|_2\\
			&\leq 2 \|P_kHP_k\| \|V_k-\hat{V}_k \|_2\leq 2\, k\,\frac{1}{k^2}
			\to 0 \quad \text{as } k\to\infty.
		\end{align}
		Thus, taking the limit in \eqref{eq:conv-ineq-8} yields
		\begin{align*}
			\|\cdot\|_1-\lim_{k\to\infty} e^{-\frac{1}{2} H} [[V_k, H], V_k] e^{-\frac{1}{2} H} = e^{-\frac{1}{2} H} [[V, H], V] e^{-\frac{1}{2} H},
		\end{align*}
		implying \eqref{eq:conv-ineq-7}. Finally, passing to the limit $k\to\infty$ in \eqref{eq:conv-ineq-6} and applying \eqref{eq:conv-ineq-7} along with the pointwise convergence $g_k\rightarrow f$ establishes \eqref{eq:conv-ineq-2}.
		
		\medskip

		{\it Case 4:} $H\geq mI$, $V$ is symmetric and $H$-bounded, $[H,V]\in\bh$ and $V(H-i)^{-1}\in\Scal^1$.
		
		\noindent	
		The convexity of $f$, along with \eqref{eq:conv-ineq-1}, was established in the course of the proof of {\it Case 3}.
		
		Let $P_k$, $Q_k$, $V_k$, $\hat V_k$ be as in {\it Case 3,} with $Q_k$ chosen by Lemma~\ref{aprlem1} with $n=1$ and $F=\{1\}$ so that $\|V_k-\hat V_k\|_1<1/k^2$. By \eqref{eq-comm-bdd}, $[H, V_k] \in \bh$. By spectral analysis analogous to the one in {\it Case 3}, $[H, \hat{V}_k] \to [H,V]$ strongly, $\|V_k(H-i)^{-1}- V(H-i)^{-1}\|_1\to 0$, $\|(H-i)^{-1}V_k-(H-i)^{-1}V\|_1\to 0$, and $\|[V_k, H] - [\hat{V}_k, H]\|_1 \to 0$.
		An argument similar to that in the proof of \eqref{eq:conv-ineq-2} in {\it Case 3} completes the proof of \eqref{eq:conv-ineq-2} in {\it Case 4}.
	\end{proof}

	In the sequel, we will need the following lemma, which is a standard consequence of Duhamel's formula for operator exponentials.
	
	\begin{lma}\label{lem:duham}
		Let $H, V$ be self-adjoint operators in $\Hcal$ such that $H\geq mI$ for some $m \in\R$ and $V \in \Scal^1$. Then, for every $t\in[-1,1]$,
		\begin{align*}
			\frac{d}{ds} e^{-(H+sV)}\Big|_{s=t} = -\int_{0}^{1} e^{-s(H+tV)} V e^{-(1-s)(H+tV)} \, ds,
		\end{align*}
		where the derivative is evaluated in the $\|\cdot\|_1$ topology.
	\end{lma}
	
	\begin{lma}\label{lem:corr-2}
		Let $m\in\R$, and let $H, V$ be self-adjoint operators in $\Hcal$ satisfying $H\geq mI$ with $V\in\Scal^1$. Define $Z:[-1,1]\to \C$ by
		\begin{align}\label{eqq-1}
			Z(t) = \Tr\bigl[e^{-(H + tV)} - e^{-H}\bigr].
		\end{align}
		Then $Z''$ is convex on $[-1,1]$ and satisfies
		\begin{align}\label{eq-2}
			Z''(0) \leq \frac12\left(\Tr\left(V e^{-H + V}\right)-\Tr\left(V e^{-H - V}\right)\right).
		\end{align}
	\end{lma}
	\begin{proof}
		By Theorem \ref{thm:Derivative-sign-bdd}, $Z^{(4)}(t)\geq 0$ for all $t\in [-1,1]$. This implies that $Z''$ is convex on $[-1,1]$. By Lemma \ref{lem:duham} and the cyclicity of the trace, $Z'(1)= -\Tr\left(V e^{-H - V}\right)$, $Z'(-1)= -\Tr\left(V e^{-H + V}\right)$. By the convexity of $Z''$,
\begin{align*}
Z''(0) \leq \frac{1}{2} \int_{-1}^{1} Z''(s)\, ds = \frac{1}{2} \big(Z'(1) - Z'(-1)\big),
\end{align*}
which completes the proof.
\end{proof}
	
	The following result is an analog of \cite[Theorem 7.1]{LNR} and \cite[Theorem~2]{DNN} for operators with essential spectrum.
	
	\begin{thm}\label{thm:corr-3}
Let $H$ be self-adjoint in $\Hcal$ such that $H\geq mI$ for some $m\in\R$. Let $V$ be symmetric and relatively $H$-bounded such that $[H, V]\in\bh$, and $V(H-i)^{-1}\in\Scal^1$. Then,
		\begin{align}\label{eq-3}
			\Tr\bigl(e^{-\frac{1}{2}H} V^2 e^{-\frac{1}{2}H}\bigr) \leq  \frac12\left(\Tr\left(V e^{-H + V}\right)-\Tr\left(V e^{-H - V}\right)\right) + \frac{1}{4} \Tr\big(e^{-\frac{1}{2} H} [[V, H], V] e^{-\frac{1}{2} H}\big).
		\end{align}
	\end{thm}
	\begin{proof}
		First we prove the result in the case $V\in\Scal^1$ and then extend it to the general case.
		
		\smallskip
		
		{\it Case 1:} $[H, V]\in\bh$ and $V\in\Scal^1$.
		
		\noindent
		By Theorem \ref{thm:corr-1}(ii),
		\begin{align}\label{eq-4}
			\Tr\bigl(e^{-\frac{1}{2}H} V^2 e^{-\frac{1}{2}H} \bigr) \leq \int_{0}^{1}\Tr(V e^{-sH} V e^{-(1-s)H})\, ds + \frac{1}{4} \Tr\big(e^{-\frac{1}{2} H} [[V, H], V] e^{-\frac{1}{2} H}\big).
		\end{align}
		Let $Z$ be given by \eqref{eqq-1}.	By Lemma \ref{lem:duham},
		\begin{align*}
			Z'(t)=&\, \Tr\big(- V e^{-(H+tV)}\big),\\
			Z''(t)=&\, \Tr\big( -V \frac{d}{ds}e^{-(H+sV)}\big|_{s=t}\big)
			=\int_{0}^{1} \Tr\big(V e^{-s(H+tV)} V e^{-(1-s)(H+tV)}\big) \, ds.
		\end{align*}
		Applying \eqref{eq-2} (see Lemma \ref{lem:corr-2}) in \eqref{eq-4} completes the proof of \eqref{eq-3} under the assumptions of {\it Case 1}.
		
		\smallskip
		
		{\it Case 2:} $[H, V]\in\bh$ and $V(H-i)^{-1}\in\Scal^1$.

		\noindent
		Let $\{V_k\}$ be a sequence of finite-rank self-adjoint operators provided by Lemma \ref{aprlem1}. By \eqref{eq-comm-bdd}, $[H,V_k]\in\bh$. Hence, applying \eqref{eq-3} to the pair $(H,V_k)$ gives
\begin{align}\label{eq-5}
\Tr\bigl(e^{-\frac{1}{2}H} V_k^2 e^{-\frac{1}{2}H}\bigr) \leq
\frac12\left(\Tr\left(V_k e^{-H + V_k}\right)- \Tr\left(V_k e^{-H - V_k}\right)\right)
+ \frac{1}{4} \Tr\big(e^{-\frac{1}{2} H} [[V_k, H], V_k] e^{-\frac{1}{2} H}\big).
		\end{align}
		By the same argument as in {\it Case 4} of the proof of Theorem \ref{thm:corr-1}, the sequences $\{e^{-\frac{1}{2}H} V_k^2 e^{-\frac{1}{2}H}\}$ and $\{e^{-\frac{1}{2} H} [[V_k, H], V_k] e^{-\frac{1}{2} H}\}$ converge in $\|\cdot\|_1$ to $ e^{-\frac{1}{2}H} V^2 e^{-\frac{1}{2}H}$ and $e^{-\frac{1}{2} H} [[V, H], V] e^{-\frac{1}{2} H}$, respectively, as $k\to\infty$.
		Replacing the Schatten norm with the operator norm in \eqref{eq-new-0} implies the strong resolvent convergence of $\{H\pm V_k\}$ to $H\pm V$. By Proposition \ref{prop-comp-lwbd}\eqref{rel-comp}, the $H$-bound of $V$ and $V_k$ is $0$ for every $k\in\N$.
Since $H\ge mI$, by Proposition \ref{prop-comp-lwbd}\eqref{lower-bd},
	$$H\pm W\geq \Big(m - \max\big\{ a\vert{}m\vert{} + b, b(1-a)^{-1} \big\}\Big) I$$
	for any $W \in \{V, V_k\}$, $k\in\N$, where $a \in [0, 1)$ and $b \in [0, \infty)$ satisfy \eqref{Hbound}. Hence,
		\begin{align*}
			\sigma(H+V_k)\cup\sigma(H-V_k) \cup \sigma(H+V) \cup \sigma(H-V)\subset \big[ m - \max\big\{ a|m| + b, b(1-a)^{-1} \big\}, \infty\big),
		\end{align*}
for every $k\in\N.$ Consequently, by \cite[Theorem VIII.20(b)]{ReedSimonI}, $(H\pm V_k-i)e^{-H\mp V_k}$ converges strongly to $(H\pm V-i)e^{-H\mp V}$.
		Since by \eqref{eq-new-1} the sequence $\{V_k(H\pm V_k-i)^{-1}\}$ converges to $V(H\pm V-i)^{-1}$ in $\|\cdot\|_1$,
		\begin{align*}
			V_k e^{-H\mp V_k} &= V_k(H\pm V_k-i)^{-1}(H\pm V_k-i)e^{-H\mp V_k}\\
			&\to V(H\pm V-i)^{-1}(H\pm V-i)e^{-H\mp V} = V e^{-H\mp V} \quad \text{ in } \|\cdot\|_1 \text{ as } k \to \infty.
		\end{align*}
		Passing to the limit as $k\to\infty$ in \eqref{eq-5} yields the desired inequality \eqref{eq-3}.
	\end{proof}
	\medskip
	
	An explicit formula for $\Tr\left(e^{-H} [[V, H], V]\right)$ in terms of the eigenvalues and eigenvectors of $H$ was derived in \cite[Eq.~(7.21)]{DNN}, provided $e^{-H}\in\Scal^1$ and $V\in\bh$. The next result establishes a representation and an estimate for $\Tr\left(e^{-H} [[V, H], V]\right)$ in terms of the spectral data of an operator $H$ with non-empty essential spectrum.

	\begin{lma}\label{Thm55}
	Let $H$ be self-adjoint in $\Hcal$ such that $H\geq mI$ for some $m\in\R$. Let $V$ be symmetric and relatively $H$-bounded such that $[H, V] \in \Scal^2$.
			\begin{enumerate}[(i)]
			\item If $V\in \Scal^2$, then
			\begin{align*}
				\Tr\left(e^{-H} [[V, H], V]\right)
				= \int_{\R} \int_{\R} (y - x)\big(e^{-x} - e^{-y}\big)\, \Tr\big(E_H(dx)\, V\, E_H(dy)\, V\big)
			\end{align*}
			and
			\begin{align}\label{eq-7}
				\Tr\left(e^{-\frac{1}{2}H} [[V, H], V] e^{-\frac{1}{2}H}\right) \leq e^{-m} \, \|[H, V]\|_2^2.
			\end{align}
			
			\item If $ V(H - i)^{-1}\in \Scal^2$, then the estimate \eqref{eq-7} holds.
		\end{enumerate}
	\end{lma}
	
	\begin{proof}
		
		{\it Case 1:} $V,\, [H, V]\in \Scal^2$.
		
		\noindent
		By \cite[Theorem 4.1, Lemma 4.7]{ds}, the set function $\nu$ defined by $$\nu(A_1\times A_2)=\Tr(E_H(A_1)VE_H(A_2)V)$$ for all Borel subsets $A_1,A_2$ of $\R$ extends to a nonnegative Borel measure on $\R^2$. Therefore, for every Borel function $\phi:\R^2\to[0,\infty)$, we can consider the scalar integral
		$$\int_\R \int_\R\phi(x,y)\,d\nu(x,y),$$ allowing the value $\infty$.
		We have
		\begin{align}\label{DOIlim}
			\lim_{n\rightarrow\infty}\int_{[-n,n]}\int_{[-n,n]}\phi(x,y)\,d\nu(x,y)
			=\int_\R \int_\R\phi(x,y)\,d\nu(x,y).
		\end{align}
		Let $\phi:\R^2\to [0,\infty)$ be given by
		\[\phi(x,y):=(y-x)(e^{-x}-e^{-y})=(y-x)^2\int_0^1e^{-((1-s)x+sy)}\,ds.\]
		Let $P_n=E_H([-n,n])$, $n\in\N$, $H_n=HP_n$, $V_n=P_nVP_n$, $\nu_n=\nu|_{[-n,n]\times[-n,n]}$. Then, by \cite[Proposition 3.2.3]{ST19},
		\begin{align}
			\label{Arep}	
			&\Tr(e^{-H_n}[[V_n,H_n],V_n])\\
			\nonumber
			&=2\Tr(e^{-H_n}V_nH_nV_n)-\Tr(e^{-H_n}H_nV_n^2)-\Tr(e^{-H_n}V_n^2H_n)\\
			\nonumber
			&=\int_\R\int_\R e^{-x}y\,d\nu_n(x,y)+\int_\R\int_\R xe^{-y}\,d\nu_n(x,y)
			-\int_\R\int_\R xe^{-x}\,d\nu_n(x,y)-\int_\R\int_\R ye^{-y}\,d\nu_n(x,y)\\
			\nonumber
			&=\int_\R\int_\R\phi(x,y)\,d\nu_n(x,y)=\int_{[-n,n]}\int_{[-n,n]}\phi(x,y)\,d\nu(x,y).
		\end{align}
		Since $V,\, [H, V]\in \Scal^2$,
		\[[V_n,H_n]=P_n[V,H]P_n\rightarrow [V,H] \text{ in } \Scal^2.\]
		Since $H\ge mI$,
		\begin{align}
			\label{rep3}
			\lim_{n\rightarrow\infty}\Tr(e^{-H_n}[[V_n,H_n],V_n])=\Tr(e^{-H}[[V,H],V]).
		\end{align}
		Combining \eqref{DOIlim}, \eqref{Arep}, and \eqref{rep3} yields
		\begin{align*}
			\Tr(e^{-H}[[V,H],V])=\int_\R \int_\R\phi(x,y)\,d\nu(x,y).
		\end{align*}
		Note that $\nu$ is supported on $[m,\infty)^2$ (since $H\ge mI$) and $0\le\phi(x,y)\le e^{-m}(x-y)^2$ for all $(x,y)\in [m,\infty)^2$. Hence,
		\begin{align*}
			0\le\Tr(e^{-H}[[V,H],V])&\le e^{-m}\int_\R \int_\R (x-y)^2\,d\nu(x,y) = \lim_{n\to\infty} e^{-m}\int_{[-n,n]} \int_{[-n,n]} (x-y)^2\,d\nu_n(x,y)\\
			&= \lim_{n\to\infty} e^{-m} \|[V_n,H_n]\|_2^2 =e^{-m}\|[V,H]\|_2^2.
		\end{align*}
		
		{\it Case 2:} $[H, V]\in\Scal^2$ and $V(H-i)^{-1}\in \Scal^2$.
		
		\noindent
		Let $\{V_k\}$ be a sequence of finite-rank self-adjoint operators given by Lemma \ref{aprlem1}. By \eqref{eq-comm-bdd} and \eqref{eq:conv-ineq-7}, $[V_k, H]$ is finite-rank and $\Tr\big(e^{-\frac{1}{2} H} [[V_k, H], V_k] e^{-\frac{1}{2} H}\big)\to \Tr\big(e^{-\frac{1}{2} H} [[V, H], V] e^{-\frac{1}{2} H}\big)$ as $k\to\infty$.
		By \eqref{eq-new-2} and \eqref{eq-new-3}, we also have $[V_k, H]\to [V, H]$ in $\|\cdotp\|_2$.
		By \eqref{eq-7} applied to the pair $(H,V_k)$,
		\begin{align*}
			0\le\Tr\big(e^{-\frac{1}{2} H} [[V_k, H], V_k] e^{-\frac{1}{2} H}\big) \le e^{-m}\|[V_k,H]\|_2^2.
		\end{align*}
		Passing to the limit as $k\to\infty$ completes the proof of \eqref{eq-7} for the pair $(H,V)$.
	\end{proof}

	By combining the results of Theorem~\ref{thm:corr-3} and Lemma~\ref{Thm55} we obtain the following essential-spectrum analog of \cite[Theorem~12]{DNN}.
	
	\begin{thm}\label{Thm66}
Let $H$ be self-adjoint in $\Hcal$ such that $H\geq mI$ for some $m\in\R$. Let $V$ be symmetric and relatively $H$-bounded such that $V(H-i)^{-1}\in\Scal^1$, and $[H, V] \in \Scal^2$. Then,
		\begin{align*}
			\Tr\bigl(e^{-\frac{1}{2}H} V^2  e^{-\frac{1}{2}H}\bigr) \leq   \frac12\left(\Tr\left(V e^{-H + V}\right)-\Tr\left(V e^{-H - V}\right)\right) + \frac{1}{4} e^{-m}\|[V,H]\|_2^2.
		\end{align*}
	\end{thm}

	\section{Examples}\label{sec:exs}
	
	In this section, we present various classes of self-adjoint operators $H$ and $V$ satisfying the hypotheses of our main results.
	
	Let $d,n\in\N$, and let $F$ be either $\{1,\ldots, n\}$ or $\{1\}$.
	
	\subsection{\it $H$ is bounded below self-adjoint, $V$ is symmetric and relatively $H$-bounded, $V(H-i)^{-p} \in \Scal^{n/p}$ for all $p\in F$.}
	
	These hypotheses on $H$ and $V$ are fundamental to Theorems \ref{thm:all_oreder_positivity-ssf}, \ref{thm:nth_order_positivity-ssf}, \ref{thm:BMV-unbdd}, and \ref{bmvsn}, and, if $n=1$, to Theorem \ref{thm2:BMV-unbdd}. Below we provide illustrative examples of $H$ and $V$ satisfying the aforementioned properties. By Lemma~\ref{lem:KR}\eqref{item:properties rel bdd: Kato-Rellich}, the respective $H+tV$ is self-adjoint on $\dom(H+tV)=\dom(H)$ for every $t\in\R$.
	
	\medskip
	
	$\blacktriangleright$ ~{\it Laplace operator on Euclidean spaces.}
	
	\noindent
	A canonical choice for the unperturbed Hamiltonian $H$ is the negative Laplacian $H = -\Delta$ defined in $L^2(\R^d)$, where $$\Delta := \sum_{k=1}^d \frac{\partial^2}{\partial x_k^2},$$ which gives rise to standard Schr\"odinger operators.
	A standard choice for the perturbation $V$ is the multiplication operator $V = M_f$ associated with a measurable real-valued function $f$ on $\R^d$.  The operator $M_f$ is self-adjoint on its maximal domain $$\dom(M_f)=\big\{h\in L^2(\R^d):\;f h\in L^2(\R^d)\big\}.$$
	
	We will use the space $\ell^p(L^2)(\R^d)$, $p \in [1,2]$, defined by
	$$\ell^p(L^2)(\R^d):=\Big\{ f : \R^d \to \C \text{ measurable }: \sum_{k\in\mathbb{Z}^d} \Big(\|f \upharpoonright_{(0,1)^d+k}\|_{L^2((0,1)^d+k)}^p\Big)<\infty\Big\}.$$
	and the space
	\begin{align*}
		H^{r}(\R^d):=\{f \in L^2(\R^d) : \partial^{\alpha}\!f \in  L^2(\R^d) \text{ for } 0 < |\alpha| \leq r\},
	\end{align*}
	where $r\in\N$ and $\partial^\alpha = \partial_1^{\alpha_1}\cdots\partial_d^{\alpha_d}$ denotes the distributional derivative of order $|\alpha|$.

(a)	Let $d, n \in \N$ with $n > d/2$, and let $f$ be a real-valued function such that
	\begin{align*}
		f\in\begin{cases}
			\ell^1(L^2)(\R^d) &\text{if }  d=1,\, n=1,\\
			L^2(\R^d) \cap L^n(\R^d) &\text{if }  1\le d/2< 2 \leq n < \infty,\\
			L^n(\R^d) &\text{if } 2 \le d/2 < n < \infty.\end{cases}
	\end{align*}
	By \cite[Theorem XIII.96]{ReedSimonIV}, $M_f$ is relatively $-\Delta$-bounded and the relative $-\Delta$-bound of $M_{f}$ is $0$. By \cite[Theorem III.3]{S21}, we have
	\begin{align*}
		M_{f}(-\Delta - i)^{-1} \in \Scal^{n}.
	\end{align*}
 Therefore, the pair $(H = -\Delta, V = M_{f})$ satisfies the hypotheses of Theorem \ref{thm:nth_order_positivity-ssf}. If in addition, $f\geq 0$, then the pair $(H=-\Delta, V=M_{f})$ satisfies the hypotheses of Theorem \ref{bmvsn}, and in case $n=d=1$, it also satisfies the hypothesis of Theorem \ref{thm2:BMV-unbdd}.
	
	\medskip
	
(b)	Let $d, n \in \N$ with $n > d/2$, and suppose that $f\in \ell^1(L^2)(\R^d)\cap \big(L^\infty(\R^d) \cup L^n(\R^d)\big)$ is real-valued. If $f$ is bounded, then  by \cite[Theorem III.4]{S21}, we have
	\begin{align}\label{i2}
		M_{f}(-\Delta - i)^{-p} \in \Scal^{n/p} \text{ for all } p \in F.
	\end{align}
	If $f$ is unbounded, then $ f\in \ell^1(L^2)(\R^d) \cap L^n(\R^d)$ implies that
	$$\{\|f\upharpoonright_{(0,1)^d+k}\|_{L^2((0,1)^d+k)}\}_{k\in\Z^d}\in \ell^1(\Z^d)\subset \bigcap\limits_{p\in [1,2 ]} \ell^p(\Z^d),$$
	consequently $f \in \bigcap\limits_{p\in [1,2 ]} \ell^p(L^2)(\R^d)$. In particular, $f\in  \ell^2(L^2)(\R^d)$, so $f\in L^2(\R^d)$. Consequently, by $L^p$ interpolation based on H\"older's inequality, $$f\in L^2(\R^d)\cap L^n(\R^d)\subset \bigcap\limits_{p\in [2\wedge n, 2\vee n]} L^p(\R^d).$$
	Therefore,
	\begin{align}\label{lpq}
		f\in \ell^1(L^2)(\R^d) \cap L^n(\R^d) \subset \Big(\bigcap\limits_{p\in [1,2 ]} \ell^p(L^2)(\R^d)\Big)\, \bigcap\, \Big(\bigcap\limits_{p\in [2\wedge n, 2\vee n]} L^p(\R^d)\Big).
	\end{align}
By the case $f\in L^2(\R^d)\cap L^n(\R^d)$ discussed above, $M_f$ is relatively $-\Delta$-bounded and the $-\Delta$-bound of $M_{f}$ is equal to $0$.
	
	For each $p\in F$, define $h_p(x)=(\|x\|^2-i)^{-p}$ on $\R^d$.	Note that
	\begin{align*}
		h_p\in
		\begin{cases}
			L^{n/p}(\R^d) & \text{ for } p=1,\ldots, \lfloor n/2 \rfloor,\\
			\ell^{n/p}(L^2)(\R^d) & \text{ for } p= \lfloor n/2 \rfloor+ 1,\ldots, n.
		\end{cases}
	\end{align*}
	Consequently, by Birman-Solomyak’s inequality \cite[Theorem 4.5]{Simon-tr-id} and Kato-Seiler-Simon's inequality \cite[Theorem 4.1]{Simon-tr-id},
	\begin{align*}
		M_{f}(-\Delta-i)^{-p}=f(x)h_p(-i\nabla) \in \Scal^{n/p}\quad \text{for all } p\in F,
	\end{align*}
	implying \eqref{i2}.

	Therefore, the pair $(H=-\Delta, V=M_{f})$ satisfies the hypotheses of Theorem \ref{thm:all_oreder_positivity-ssf}. If, in addition, $f\geq 0$, then the pair $(H=-\Delta, V=M_{f})$ satisfies the hypotheses of Theorems \ref{thm:BMV-unbdd}, \ref{bmvsn},
	and, in case $n=d=1$, it also satisfies the hypothesis of Theorem \ref{thm2:BMV-unbdd}.
	
	\medskip
	
	$\blacktriangleright$ ~{\it Fractional Laplacian operators on Euclidean spaces.}
	
	\noindent
	Here we use the notation introduced in the previous example.
	
	\smallskip
	
	\noindent
	Let $d, n \in\N$ with $d<2n$, let $s \in \left(\frac{d}{2n}, 1\right)$. For a given $p \in F$, consider the fractional Laplacian $H = (-\Delta)^s$. Let $f \in \ell^1(L^2)(\R^d)\cap L^\infty(\R^d)$ be real-valued, let $V = M_{f}$, and define $h_p(x) = (\|x\|^{2s} - i)^{-p}$ on $\R^d$. Then, by Birman-Solomyak’s inequality \cite[Theorem 4.5]{Simon-tr-id} and Kato-Seiler-Simon's inequality \cite[Theorem 4.1]{Simon-tr-id},
	\begin{align*}
		V(H-i)^{-p}=f(x)h_p(-i\nabla) \in \Scal^{n/p} \quad \text{for all } p \in F.
	\end{align*}
	Therefore, the pair $(H, V)$ satisfies the hypotheses of Theorems \ref{thm:all_oreder_positivity-ssf}, \ref{thm:nth_order_positivity-ssf}. If in addition, $f\geq 0$, then the pair $(H, V)$ satisfies the hypotheses of Theorems \ref{thm:BMV-unbdd} and \ref{bmvsn} and, in case $n=d=1$, it also satisfies the hypothesis of Theorem \ref{thm2:BMV-unbdd}.
	
	\medskip

	$\blacktriangleright$ ~{\it Magnetic anharmonic oscillator.}
	
	\noindent
	Let $d,m,n\in\N$, $\kappa_{d,m},\theta\in\R$ satisfy
	\begin{equation*}
		\kappa_{d,m}:=\frac{d(m+1)}{2m}<n,\quad 0<\theta<\min\left\{\frac12,\,1-\frac{\kappa_{d,m}}{n}\right\}.
	\end{equation*}
	Let $A=(A_1,\ldots,A_d)\in C_c^\infty(\R^d;\R^d)$, the space of all compactly supported smooth functions from $\R^d$ to $\R^d$, and $w_m(x)=1+\|x\|^{2m}$.
	By \cite[Proposition~10.10]{Haslinger2023}, the symmetric operator
	\[T_{A,m}:=\sum_{j=1}^d(-i\partial_{x_j}-M_{A_j})^2+M_{w_m},
	\quad \dom(T_{A,m})=C_c^\infty(\R^d),\]
	is essentially self-adjoint. Let
	\begin{align}\label{dfn-mgnharop}
		H_{A,m}:=\overline{T_{A,m}}
	\end{align} be the (operator) closure of $T_{A,m}$.
	Since $A\in L^\infty(\R^d;\R^d)$, the vector operator $-i\nabla-A$, with the derivatives understood in the distributional sense, is closed on its maximal domain $ H^1(\R^d)$.
	The multiplication operator $M_{w_m}^{1/2}$ is closed on its maximal domain $L^2(\R^d,w_m(x)\,dx)$.
	Consider the quadratic form \[\mathfrak h_{A,m}[h]=\|(-i\nabla-A)h\|^2+\|w_m^{1/2}h\|^2,
	\quad \dom(\mathfrak h_{A,m})=H^1(\R^d)\cap L^2(\R^d,w_m(x)\,dx).\]
	By the closedness of $-i\nabla-A$ and $M_{w_m}^{1/2}$, the column operator $\mathcal{C}$ given by $\mathcal{C}h:=((-i\nabla-A)h,M_{w_m}^{1/2}h)$ is closed on $\dom(\mathfrak h_{A,m})$.
	Hence, $\dom(\mathfrak h_{A,m})$ is complete in the graph norm of $\mathcal{C}$, which is also the form norm of $\mathfrak h_{A,m}$ given by $(\mathfrak h_{A,m}[h]+\|h\|^2)^{1/2}$. By \cite[Chapter~VI, Section~1.3, Theorem~1.11]{Kato1995}, $\mathfrak h_{A,m}$ is closed. By \cite[Theorem~9.3]{Simon2018Kato}, $C_c^\infty(\R^d)$ is a form core for $\mathfrak h_{A,m}$. Since integration by parts gives $\mathfrak h_{A,m}[h]=\left<T_{A,m}h,h\right>$ for $h\in C_c^\infty(\R^d)$ (see, e.g., the proof of \cite[Proposition~10.7]{Haslinger2023}), the quadratic form of $T_{A,m}$ is $\mathfrak h_{A,m}|_{C_c^\infty(\R^d)}$. Since a form core is, by definition \cite[Chapter~VI, Section~1.4, p.~317]{Kato1995}, dense in the form domain with respect to the form norm, the closure of this restriction is $\mathfrak h_{A,m}$. Since $\dom(\mathfrak h_{A,m})$ contains $C_c^\infty(\R^d)$, it is dense in $L^2(\R^d)$. Since $\mathfrak h_{A,m}$ is densely defined, closed, symmetric,
	and nonnegative, Kato's first representation theorem
	\cite[Chapter~VI, Section~2.1, Theorem~2.1]{Kato1995} along with \cite[Chapter~VI, Section~2.2, Theorem~2.6]{Kato1995}
	determines a unique nonnegative self-adjoint operator associated
	with $\mathfrak h_{A,m}$. This operator is the Friedrichs extension
	of $T_{A,m}$ (see \cite[Theorem~X.23]{ReedSimonII} or \cite[Chapter~VI, Section~2.3]{Kato1995}) and hence, by essential self-adjointness, equals
	$H_{A,m}=\overline{T_{A,m}}$. Kato's second representation theorem
	\cite[Chapter~VI, Section~2.6, Theorem~2.23]{Kato1995}
	then gives $\dom(H_{A,m}^{1/2})=\dom(\mathfrak h_{A,m})$.
	Since $\left<H_{A,m}h,h\right>\ge\int_{\R^d}w_m(x)|h(x)|^2\,dx$ for every $h\in C_c^\infty(\R^d)$ and $w_m(x)\to\infty$ as $x\to\infty$,
	\cite[Theorem~11.9(d)]{Haslinger2023} implies that $H_{A,m}$ has
	compact resolvent.
	
	Let
	\[V=M_{w_m^\theta}=M_{w_m}^\theta.\]
	Below we justify that $V$ is relatively $H_{A,m}$-bounded and that $V(H_{A,m}-i)^{-p}\in\S^{n/p}$, $p=1,\ldots,n$.
	
	Let
	\begin{align*}
		&T_m=-\Delta+M_{w_m},\quad \dom(T_m)=H^2(\R^d)\cap\dom(M_{w_m}),\\
		&\mathfrak t_m[h]=\|\nabla h\|^2+\|w_m^{1/2}h\|^2, \quad h\in\dom(\mathfrak h_{A,m}).
	\end{align*}
	By \cite[Theorem~4.1]{Dung}, $T_m$ is self-adjoint on this domain.
	Since $-i\nabla h=(-i\nabla-A)h+Ah$ for $h\in \dom(\mathfrak h_{A,m})$ and $w_m\ge1$, we obtain
	\begin{align}
		\label{form_compare}  \nonumber
		\mathfrak t_m[h]&\le 2\|(-i\nabla-A)h\|^2+2\|A\|_\infty^2\|h\|^2 + \|w_m^{1/2}h\|^2\\
		\nn&\le 2\|(-i\nabla-A)h\|^2+(2\|A\|_\infty^2 + 1)\|w_m^{1/2}h\|^2\\
		&\le C_A\,\mathfrak h_{A,m}[h],
		\quad\text{where}\quad C_A:=2+2\|A\|_\infty^2.
	\end{align}
	By \cite[Corollary~5.4(b)]{ChDeRu21},
	$\bigl((-\Delta)^\ell+M_{w_k}\bigr)^{-\mu}\in\S^r$ when $\mu>\frac{d(k+\ell)}{2k\ell r}$.
	Taking $k=m$, $\ell=1$, $\mu=s$, and $r=q$ gives
	\[T_m^{-s}\in\S^q \quad\text{if}\quad sq>\frac{d(m+1)}{2m}=\kappa_{d,m}.\]
	Let $\{\mu_j\}_j$ and $\{\lambda_j\}_j$ denote the nondecreasing sequences of eigenvalues of $T_m$ and $H_{A,m}$, respectively, with each eigenvalue repeated according to its multiplicity.
	By \eqref{form_compare} and the min-max principle \cite[Theorem~XIII.1]{ReedSimonIV},
	$\mu_j\le C_A\lambda_j$.  Hence,
	$\sum_j\lambda_j^{-sq} \le C_A^{sq}\sum_j\mu_j^{-sq}<\infty$, so
	\begin{equation}\label{HAmSq}
		H_{A,m}^{-s}\in\S^q \quad\text{whenever}\quad sq>\kappa_{d,m}.
	\end{equation}
	Since $\|M_{w_m}^{1/2}h\|^2\leq\|H_{A,m}^{1/2}h\|^2$ for each
	$h\in\dom(H_{A,m}^{1/2})\subset \dom(M_{w_m}^{1/2})$ and $0<2\theta<1$, the Heinz inequality \cite[Proposition~10.14]{Schmuedgen2012} gives $\|M_{w_m}^{\theta}h\|^2\leq\|H_{A,m}^{\theta}h\|^2$ for every $h\in\dom(H_{A,m}^{\theta})$. That is,
	\begin{equation}\label{VHtheta}
		\|Vh\|\le\|H_{A,m}^\theta h\|,\quad h\in\dom(H_{A,m}^{\theta}).
	\end{equation}
	Since the function $t\mapsto t^\theta-\varepsilon t$ attains its maximum on $[0,\infty)$ at $\big(\frac{\theta}{\varepsilon}\big)^{\frac{1}{1-\theta}}$, we have that
	for every $\varepsilon>0$ and $t\ge1$, $t^\theta\le\varepsilon t+C_{\varepsilon,\theta}$,
	where
	$C_{\varepsilon,\theta}=\big(\frac{\theta}{\varepsilon}\big)^{\frac{\theta}{1-\theta}}
	-\varepsilon\big(\frac{\theta}{\varepsilon}\big)^{\frac{1}{1-\theta}}
	=(1-\theta)(\frac{\theta}{\varepsilon})^{\frac{\theta}{1-\theta}}$.
	Combining the latter with \eqref{VHtheta} gives
	\[\|Vh\|\le\varepsilon\|H_{A,m}h\|+C_{\varepsilon,\theta}\|h\|,
	\quad h\in\dom(H_{A,m}).\]
	Thus $V$ has relative $H_{A,m}$-bound $0$. Note also that $VH_{A,m}^{-\theta}$ is bounded.
	
	Since $(p-\theta)n/p=n\left(1-\theta/p\right)\ge n(1-\theta)>\kappa_{d,m}$,
	applying \eqref{HAmSq} with $s=p-\theta$ and $q=n/p$ gives
	\[H_{A,m}^{-(p-\theta)}\in\S^{n/p}.\]
	Let $\{\sigma_j\}_j$ denote the nonincreasing sequence of singular values of $H_{A,m}^\theta(H_{A,m}-i)^{-p}$. Since
	$\frac{t^\theta}{|t-i|^p}\le t^{-(p-\theta)}$ for $t\ge 1$, we have
	$\sigma_j\le \lambda_j^{-(p-\theta)}$, and hence,
	\[H_{A,m}^\theta(H_{A,m}-i)^{-p}\in\S^{n/p}.\]
	Therefore,
	\[V(H_{A,m}-i)^{-p}=\big(VH_{A,m}^{-\theta}\big)\big(H_{A,m}^\theta(H_{A,m}-i)^{-p}\big)
	\in\S^{n/p},\quad p=1,\ldots,n.\]

	Therefore, the pair $(H=H_{A,m},V)$ satisfies the hypotheses of Theorems \ref{thm:all_oreder_positivity-ssf}, \ref{thm:BMV-unbdd}, \ref{bmvsn}, and, in case $n=d=1$, it also satisfies the hypothesis of Theorem \ref{thm2:BMV-unbdd}.

	\subsection{\it $V\in\bh$, $[H,V]\in\Scal^2$ and $V(H-i)^{-1}\in \Scal^1$.}
		
	These hypotheses on $H$ and $V$ are fundamental to  Theorems \ref{thm:corr-3} and \ref{Thm66}. We  provide illustrative examples below.
	
	\medskip
	
	$\blacktriangleright$ ~{\it Lattice Laplacians.}
	
	\noindent
	Let $\H=\ell^2(\Z^d)$.
	Let $H$ be the nearest-neighbor operator given by
	$$(\Delta_d\,x)_{\vec{k}}=\sum\limits_{\vec{j}:\; \|\vec{j}-\vec{k}\|=1}\!x_{\vec{j}},\quad x=(x_{\vec{k}})_{\vec{k}}\in\ell^2(\Z^d),$$
	or the lattice Laplacian given by
	$$\mathsf L_{d}=2dI - \Delta_{d}.$$
	Let $(v_{\vec{k}})_{\vec{k}}$ be a real-valued sequence in $\ell^1(\Z^d)$, and let $V$ be defined by
	$$(Vx)_{\vec{k}}=v_{\vec{k}}\,x_{\vec{k}},\quad x=(x_{\vec{k}})_{\vec{k}}\in\ell^2(\Z^d).$$
	In the case $d=1$, $\Delta_1$ is the full-line Jacobi matrix acting on $\ell^2(\Z)$ by $(\Delta_1 x)_k=x_{k-1}+x_{k+1}$, $k\in\Z$, and $V={\rm diag}(v)$.
	
	One can also consider $H=\alpha I+\beta\Delta_{d}$ for some $\alpha,\beta\in\R$. Another well-studied model is the half-line Jacobi matrix $J$ on $\ell^2(\N_0)$ defined by $(Jx)_0=x_1$, $(Jx)_k=x_{k-1}+x_{k+1}$, $k\in\N$, and perturbed by $V={\rm diag}(v)$.
	A more detailed discussion of discrete Laplacians and their applications can be found in, for instance, \cite{CL90,KLW21}. The theory of more general weighted graph Laplacians can be found in \cite{KLW21}.
	
	It is easy to see that $H=H^*\in\mathcal{B}(\ell^2(\Z^d))$ and $V=V^*\in\Scal^1(\ell^2(\Z^d))$.
	Consequently, the pair $(H,V)$  satisfies the hypotheses of Theorems \ref{thm:all_oreder_positivity-ssf}, \ref{thm:nth_order_positivity-ssf}, \ref{thm:corr-3} and \ref{Thm66}.
	Furthermore, if $v_{\vec{k}}\geq 0$ for every $\vec{k}\in\Z^d$, then $(H,V)$ satisfies the hypotheses of Theorems \ref{thm:BMV-unbdd}, \ref{thm2:BMV-unbdd}, and \ref{bmvsn}. Analogous properties hold for the pair $(J,V)$.
	
	\medskip
	
	$\blacktriangleright$ ~{\it Friedrichs-type operators.}
	
	\noindent
	Below we consider examples of a Friedrichs-type multiplication-integral operator and a regular Friedrichs--Lee block operator. A more detailed discussion of Friedrichs models can be found, for instance, in \cite[Chapter~3, Section~3.2]{Ex85} and \cite{Fr65}.
	
	(a) Let $d\geq 1$ and $\Hcal=L^2(\R^d)$. Let $\omega:\R^d\to\R$ be a measurable function
	and let $H_0=M_\omega$ on its maximal multiplication-operator domain.
	Then, $H_0$ is self-adjoint. Given $g\in\dom(H_0)\setminus\{0\}$ and $\lambda\in\R\setminus\{0\}$, define $V$ on $\Hcal$ by
	$(Vf)=\lambda g\int_{\R^d}f(x)\overline{g(x)}\,dx$. Note that ${\rm rank}(V)=1$ and $\ran(V)\subset\dom(H_0)$. Set $H=H_0+V$ and note that ${\rm rank}([H,V])\le 2$.
	
	(b) Let $\mathcal{K}$ be a Hilbert space and set $\Hcal=\C\oplus\mathcal{K}$. Let $K$ be a self-adjoint operator in $\mathcal{K}$, let $\varepsilon\in\R$ and define $H_0=\begin{pmatrix} \varepsilon&0\\ 0&K \end{pmatrix}$. Then, $\dom(H_0)=\C\oplus\dom(K)$. Given $g\in\dom(K)\setminus\{0\}$, define $C_{g}:\C\to\mathcal{K}$ by $C_{g}z=zg$ and let $V=\begin{pmatrix}0&C_{g}^{*}\\ C_{g}&0 \end{pmatrix}$. Note that $\ran(V)={\rm span}\{(1,0),(0,g)\}\subset\dom(H_0)$. Set $H=H_0+V$ and note that ${\rm rank}([H,V]) \leq 2$.
	
	Consequently, if $\omega$ is bounded from below a.e. in (a) and $K$ is bounded from below in (b), the pair $(H,V)$ in both (a) and (b) satisfies the hypotheses of Theorems  \ref{thm:corr-3}, and \ref{Thm66}.
	If $\omega$ is bounded from below a.e. and $\lambda >0$, then the pair $(H,V)$ defined in (a) also satisfies the hypotheses of Theorems \ref{thm:all_oreder_positivity-ssf}, \ref{thm:nth_order_positivity-ssf}, \ref{thm:BMV-unbdd}, \ref{thm2:BMV-unbdd}, and \ref{bmvsn}.
	
	\medskip

	$\blacktriangleright$ ~{\it  Anharmonic oscillator--Laplacian direct sum.}
	
	\noindent
	Before stating our example, we discuss sufficient conditions for operators $H, V$ in a separable Hilbert space to satisfy $[V, H]\in\Scal^2$ and $V(H-i)^{-1}\in \Scal^2$
	(or $V(H-i)^{-1}\in \Scal^1$).

	\begin{lma}\label{lem:sufficient}
		Let $H_0$ and $H_1$ be two self-adjoint operators in $\Hcal$ such that $(H_0-i)^{-1}\in \Scal^2$. Let $W\in \Scal^2(\Hcal)$ and for a given $n\in\N$, define $W_n=E_{H_1}([-n, n]) W E_{H_1}([-n, n])$. Let $P$ be a finite-rank projection such that $P(\Hcal)\subset \dom(H_0)\cap \dom(H_1)$, and let $f:\R\to\C$ be a bounded measurable function. Define the operators $H$ and $V$ on $\Hcal\oplus\Hcal$ by
		\begin{align}\label{def:H}
			H &= H_0\oplus H_1,\\
			\label{def:V}V &= \begin{pmatrix}
				f(H_0) & P\\
				P & W_n
			\end{pmatrix}.
		\end{align}
		Then,
		\begin{align*}
			[V, H] &= \begin{pmatrix}
				0 & PH_1-H_0P\\
				PH_0-H_1P & [W_n, H_1]
			\end{pmatrix}\in \Scal^2(\Hcal\oplus\Hcal)\\[10pt]
			V(H-i)^{-1} &= \begin{pmatrix}
				f(H_0)(H_0-i)^{-1} & P(H_1-i)^{-1}\\
				P(H_0-i)^{-1} & W_n(H_1-i)^{-1}
			\end{pmatrix}\in \Scal^2(\Hcal\oplus\Hcal).
		\end{align*}
		Moreover, if  $(H_0-i)^{-1}, W\in \Scal^1(\Hcal)$, then $[V, H], V(H-i)^{-1}\in \Scal^1(\Hcal\oplus\Hcal)$.
	\end{lma}
	\begin{proof}
		We note that each block of $[V, H]$ and $V(H-i)^{-1}$ is in $\Scal^2(\Hcal)$, hence $[V, H], V(H-i)^{-1}\in \Scal^2(\Hcal\oplus\Hcal)$. Moreover, if $(H_0-i)^{-1}, W\in \Scal^1(\Hcal)$, then by similar arguments $[V, H], V(H-i)^{-1}\in \Scal^1(\Hcal\oplus\Hcal)$. This completes the proof.
	\end{proof}
	
	Let $p \in \{1, 2\}$ and $d, l, k \in \N$ satisfy $2klp > d(k+l)$.
	Let $$H_0 = (-\Delta)^l + M_{v_k}, \quad v_k(x)=\|x\|^{2k},$$
	be an anharmonic oscillator defined in $L^2(\R^d)$ (see, e.g., \cite{ChDeRu21}) and let $H_1=-\Delta$. By \cite[Theorem~4.1]{Dung}, $H_0$ is self-adjoint in $L^2(\R^d)$ with
	$\dom(H_0) =  H^{2l}(\R^d)\cap\dom(M_{v_k})$
	and $H_1$ is self-adjoint in $L^2(\R^d)$ with $\dom(H_1)=H^2(\R^d)$.
	
	Then, by \cite[Corollary 5.4 (b)]{ChDeRu21}, we have $(H_0 + 1)^{-1} \in \Scal^p(L^2(\R^d))$, and consequently $(H_0 -i)^{-1} \in \Scal^p(L^2(\R^d))$. Let $f:\R \to \R$ be a bounded measurable function, and for $p\in\{1, 2\}$, let $W\in \Scal^p(L^2(\R^d))$ be self-adjoint. Let $P : L^2(\R^d) \to L^2(\R^d)$ be a finite-rank orthogonal projection such that $\ran(P) \subset \dom(H_0) \cap \dom(H_1) = \dom(H_0)$. By Lemma \ref{lem:sufficient}, the operators $H$ and $V$ given by \eqref{def:H} and \eqref{def:V} satisfy
	$$[V, H], \, V(H - i)^{-1} \in \Scal^p(L^2(\R^d) \oplus L^2(\R^d))$$
	for $p \in \{1, 2\}$. Thus, the pair $(H, V)$ satisfies the hypotheses of Theorems \ref{thm:corr-3} and \ref{Thm66} when $p = 1$, and of Theorem \ref{thm:corr-1} when $p = 2$. If $f\equiv 1$, then $V$ is not compact.
	
	\bigskip

	\noindent\textit{Acknowledgment}:
	C. Pradhan gratefully acknowledges support from the Fulbright-Nehru Postdoctoral Fellowship, Award No. 3130/FNPDR/2025. A. Skripka is supported in part by Simons Foundation Grant MP-TSM-00002648.

\end{document}